\documentclass[12pt,reqno]{amsart}
\usepackage{amsmath}
\usepackage{amssymb}
\usepackage{amsthm}
\usepackage{latexsym}
\usepackage{amssymb}

\newtheorem{theorem}{Theorem}[section]
\newtheorem{lemma}{Lemma}[section]
\newtheorem{corollary}{Corollary}[section]
\newtheorem{remark}{Remark}[section]
\newtheorem{definition}{Definition}[section]
\newtheorem{proposition}{Proposition}[section]
\newtheorem{example}{Example}[section]
\newtheorem{assumption}{Assumption}[section]
\numberwithin{equation}{section}
\newcommand{\bth}{\begin{theorem}}
\newcommand{\ethe}{\end{theorem}}

\newcommand{\bre}{\begin{remark}}
\newcommand{\ere}{\end{remark}}

\newcommand{\ble}{\begin{lemma}}
\newcommand{\ele}{\end{lemma}}
\newcommand{\bde}{\begin{definition}}
\newcommand{\ede}{\end{definition}}

\newcommand{\bco}{\begin{corollary}}
\newcommand{\eco}{\end{corollary}}

\newcommand{\bpr}{\begin{proposition}}
\newcommand{\epr}{\end{proposition}}

\newcommand{\bexer}{\begin{exercise}}
\newcommand{\eexer}{\end{exercise}}
\newcommand{\breh}{\begin{hint}}
\newcommand{\ereh}{\end{hint}}

\newcommand{\halmos}{\hfill \qed}

\newcommand{\bexam}{\begin{example}}
\newcommand{\eexam}{\end{example}}

\newcommand{\pr} {{\bf Proof.}}

\newcommand{\bfi}{\begin{fig}}
\newcommand{\efi}{\end{fig}}

\newcommand{\beao}{\begin{eqnarray*}}
\newcommand{\eeao}{\end{eqnarray*}\noindent}

\newcommand{\beam}{\begin{eqnarray}}
\newcommand{\eeam}{\end{eqnarray}\noindent}
\newcommand{\E}{\mathbf{E}}
\newcommand{\PP}{\mathbf{P}}
\newcommand{\nto}{n\to\infty}

\newcommand{\xto}{x\to\infty}

\newcommand{\bF}{\overline{F}}

\newcommand{\bG}{\overline{G}}

\newcommand{\bbr}{{\mathbb R}}

\newcommand{\bbb}{{\mathbb B}}

\newcommand{\bbn}{{\mathbb N}}

\newcommand{\yto}{y\to\infty}

\newcommand{\vep}{\varepsilon}

\begin{document}
\title[Joint tail of randomly weighted sums]{Joint tail of randomly weighted sums under generalized quasi asymptotic independence}

\author[ D.G. Konstantinides, C. D. Passalidis ]{ Dimitrios G. Konstantinides, Charalampos  D. Passalidis}

\address{Dept. of Statistics and Actuarial-Financial Mathematics,
University of the Aegean,
Karlovassi, GR-83 200 Samos, Greece}
\email{konstant@aegean.gr,\;sasm23002@sas.aegean.gr.}

\date{{\small \today}}

\begin{abstract}
In this paper we revisited the classical problem of max-sum equivalence of randomly weighted sums in two dimensions. In opposite to the most papers in literature, we consider that there exists some interdependence between the primary random variables, which is achieved by a combination of a new dependence structure with some two-dimensional heavy-tailed classes of distributions. Further, we introduce a new approach in two-dimensional regular varying distributions, that in contrast to well-established multivariate regularly varying distributions, is consistent with the multivariate non-linear single big jump principle. We study some closure properties of this, and of other two-dimensional classes. Our results contain the finite-time ruin probability in a two-dimensional discrete time risk model.
\end{abstract}

\maketitle
\vspace{3mm}


\section{Introduction} \label{sec.KLP.1}

Let $X_1,\,\ldots,\,X_n$, $n\in\mathbb{N}$, be real valued  random variables, with distributions $F_1,\,\ldots,\,F_n$ respectively, and we consider their tails $\bF_i(x)=1-F_i(x)=\PP[X_i>x]$, with $\bF_i(x)>0$ for any $x \in \bbr$ and any $i=1,\,\ldots,\,n$. A well-known question for heavy-tailed distributions is the validity of the asymptotic relation
\beam \label{eq.KLP.1.1}
\PP[S_n> x] \sim \sum_{i=1}^n \bF_i(x)\,,
\eeam 
as $\xto$, where $S_n :=\sum_{i=1}^n X_i$, or furthermore the validity of 
\beam \label{eq.KLP.1.2}
\PP[S_n> x] \sim \PP\left[\bigvee_{i=1}^n S_i > x \right] \sim \PP\left[\bigvee_{i=1}^n X_i > x \right] \sim \sum_{i=1}^n \bF_i(x)\,,
\eeam 
as $\xto$, where 
\beam \label{eq.KLP.1.3}
\bigvee_{i=1}^n S_i :=\max_{1\leq k \leq n} \sum_{i=1}^k X_i\,, \qquad \bigvee_{i=1}^n X_i :=\max_{1\leq i \leq n} X_i\,, 
\eeam 

For the study of the relations \eqref{eq.KLP.1.1} or \eqref{eq.KLP.1.2} two factors play crucial role. The first factor is the distribution class to which belong the  $F_1,\,\ldots,\,F_n$ and the second one is the dependence structure among these random variables. For some papers on this topic see \cite{ko:tang:2008}, \cite{tang:2008}, \cite{yang:wang:leipus:siaulys:2011}.

This question is extended in the case of randomly weighted sums, where now we have some non-negative, non-degenerate to zero random variables $\Theta_1,\,\ldots,\,\Theta_n$, hence the relation \eqref{eq.KLP.1.2} becomes
\beam \label{eq.KLP.1.4}
\PP\left[S_n^{\Theta}> x\right] \sim \PP\left[\bigvee_{i=1}^n S_i^{\Theta} > x \right] \sim \PP\left[\bigvee_{i=1}^n \Theta_i\,X_i > x \right] \sim \sum_{i=1}^n \PP[\Theta_i\,X_i>x]\,,
\eeam 
as $\xto$, where 
\beam \label{eq.KLP.1.5}
S_n^{\Theta} :=\sum_{i=1}^n \Theta_i\,X_i\,, \quad \bigvee_{i=1}^n S_i^{\Theta} :=\max_{1\leq k \leq n} \sum_{i=1}^k \Theta_i\,X_i\,, \quad \bigvee_{i=1}^n \Theta_i\,X_i :=\max_{1\leq i \leq n}\Theta_i\, X_i\,, 
\eeam 
It is easily to observe that if the random weights  $\Theta_1,\,\ldots,\,\Theta_n$ are degenerated to unity, then we return back to relation \eqref{eq.KLP.1.2}. 

For some papers about the asymptotic behavior of the randomly weighted sums see in \cite{tang:tsitsiashvili:2003}, \cite{tang:yuan:2014}, \cite{yang:leipus:siaulys:2012} and \cite{konstantinides:leipus:passalidis:siaulys:2025}.

Recently, there were several attempts to extend this problem in two dimensional set-up. Namely, was arisen the question about the validity of 
\beam \label{eq.KLP.1.6}
\PP[S_n> x\,,\;T_m> y] \sim \sum_{i=1}^n \sum_{j=1}^m \PP[X_i>x\,,\;Y_j>y ]\,,
\eeam 
as $(x,\,y) \to (\infty,\,\infty)$, or as $x\wedge y \to \infty$, where $T_m = \sum_{j=1}^m Y_j$, (with $m\in\mathbb{N}$). We observe that this is not the only way to study relation \eqref{eq.KLP.1.1} in multidimensional set up. In fact, relation \eqref{eq.KLP.1.6} is the multivariate non-linear approximation of the single big jump principle, while the expression
\beam \label{eq.KLP.1.A.1}
\PP\left[\sum_{i=1}^{n} {\bf Z}^{(i)} \in x\,A \right]\sim \sum_{i=1}^{n}\PP\left[ {\bf Z}^{(i)} \in x\,A \right]\,,
\eeam
as $\xto$, for some Borel set $A \in \bbr^{2} \setminus \{{\bf 0} \}$, with ${\bf Z}^{(i)}=(X^{(i)}\,,\;Y^{(i)})$, for any $i =1,\,\ldots,\,n$, represents the multivariate linear approximation of the single big jump principle. We notice here, that relation \eqref{eq.KLP.1.A.1} is satisfied by many well-established distributions, as for example when the ${\bf Z}^{(i)}$ are independent and identically distributed random variables with standard bi-variate regular variation, see definition in Section 2, or when they are independent, identically, multivariate subexponential distributed, see in \cite[Cor. 4.10]{samorodnitsky:sun:2016}, as also under slightly smaller multivariate distribution classes, but under more relaxed conditions than identically distributed and independence of $\{{\bf Z}^{(i)}\,,\;i\in \bbn \}$, see in \cite[Sec. 4.1]{konstantinides:passalidis:2024g}. However, all the multidimensional distributions, that satisfy relation \eqref{eq.KLP.1.A.1}, can not satisfy also relation \eqref{eq.KLP.1.6}. The multivariate non-linear single big jump principle focus on the joint tail, and simultaneously give some flexibility in several applications, as for example in continuous time risk models, where we can use two different counting processes, that have impact on the final asymptotic behavior of the ruin probability, see for example in \cite{chen:wang:wang:2013} and in \cite{xu:shen:wang:2025}. This impact disappears in the multivariate regular variation, see for example in \cite{konstantinides:li:2016}, \cite{yang:su:2023}. The extension of relation \eqref{eq.KLP.1.6}, in the weighted case, becomes
\beam \label{eq.KLP.1.7}
\PP\left[S_n^{\Theta}> x\,,\;T_m^{\Delta}> y\right] \sim \sum_{i=1}^n \sum_{j=1}^m \PP[\Theta_i\,X_i>x\,,\;\Delta_j\,Y_j>y ]\,,
\eeam
as $(x,\,y) \to (\infty,\,\infty)$, or as $x\wedge y \to \infty$, where 
\beao
T_m^{\Delta} = \sum_{j=1}^m \Delta_j\,Y_j\,,
\eeao 
or if are true the asymptotic equivalencies for the maximums of the sums in two dimensions and the jointly maximums, with the 
\beao
\bigvee_{j=1}^m T_j\,, \qquad \bigvee_{j=1}^m Y_j\,, \qquad \bigvee_{j=1}^m T_j^{\Delta}\,, \qquad \bigvee_{j=1}^m \Delta_j\,Y_j\,,
\eeao 
defined in a similar way as in relations \eqref{eq.KLP.1.3} and \eqref{eq.KLP.1.5}. The reason, for which we are interested in the joint behavior of the randomly weighted sums is the dependence, that can appear among these sums. Obviously, if $S_n^{\Theta}$, $T_m^{\Theta}$ are independent, then there is no reason to study \eqref{eq.KLP.1.7}, because the problem can be reduced in \eqref{eq.KLP.1.4} for each sum.

Relation \eqref{eq.KLP.1.7} was studied in several papers where was examined one of the following cases:

1) The main random variables $\{X_1,\,\ldots,\,X_n\}$  and $\{Y_1,\,\ldots,\,Y_m\}$ are independent sequences of random variables and the two sequences are also independent. The random weights  $\Theta_1,\,\ldots,\,\Theta_n,\, \Delta_1,\,\ldots,\,\Delta_m$ are arbitrarily dependent, non-negative and non-degenerate to zero random variables, and independent of the main variables $\{X_1,\,\ldots,\,X_n\}$  and $\{Y_1,\,\ldots,\,Y_m\}$. Therefore, the dependence between the two sums $S_n^{\Theta},\,T_m^{\Delta}$ comes from the random weights.

2) The random weights are still as in case 1), but the main random variables have dependence structure. Namely, the $\{X_1,\,\ldots,\,X_n\}$ have some dependence structure and the $\{Y_1,\,\ldots,\,Y_m\}$ have also some dependence structure, but the two sequences of $\{X_i\}$ and $\{Y_j\}$ are mutually independent, see for example in \cite{yang:chen:yuen:2024}.

3) The random weights are still as in case 1), the main random variables  $\{X_i\}$ and $\{Y_j\}$ represent each a sequence of independent random variables, but the pair $(X_i,\,Y_i)$ has some dependence structure for any $i=1,\,\ldots,\,n\wedge m$.

In case 3) there exist another distinction between the papers, in these where the dependence among pairs  $(X_i,\,Y_i)$  is 'weak', that means it is close to asymptotic independence, see for example \cite{li:2018b}, and in those with arbitrary dependence among the pairs  $(X_i,\,Y_i)$, see for example \cite{chen:yang:2019}, however in most of the cases there exist assumption of identically distributed random pairs, as for example in \cite{shen:du:2023} and \cite{shen:ge:fu:2020}.

In \cite{konstantinides:passalidis:2023} was introduced the generalized tail asymptotic independence (GTAI), see in Subsection 2.3 below, as an attempt to merge the 2) and 3) dependence structures, however as result was found that GTAI covers the case when each pair 
$(X_i,\,Y_i)$ follows a special form of weak dependence. Further in \cite{konstantinides:passalidis:2025}, assuming the sums restricted by the case when $m=n$, was established the asymptotic
\beam \label{eq.KLP.1.8} \notag
\PP[S_n> x\,,\;T_n> y]  &\sim& \PP\left[\bigvee_{i=1}^n S_i > x\,,\;\bigvee_{j=1}^n T_j > y \right] \sim\PP\left[\bigvee_{i=1}^n X_i > x\,,\;\bigvee_{j=1}^n Y_j > y\right]  \\[2mm]
&\sim& \sum_{i=1}^n \sum_{j=1}^n \PP[X_i>x\,,\;Y_j>y ]\,,
\eeam
as $x\wedge y \to \infty$, under GTAI dependence structure (and TAI dependence structure on each sequence) and under some two-dimensional distribution class. Furthermore, under some conditions on the random weights, that are independent of main random variables, the asymptotic relation \eqref{eq.KLP.1.8} is still true. 

In this  work we have two aims.

$1)$ to establish partially relation \eqref{eq.KLP.1.8} for some greater dependence structure, which was possible reducing the set-up from class $(\mathcal{D}\cap \mathcal{L})^{(2)}$ to class $\mathcal{C}^{(2)}$. Further in our results  it is not necessary $n=m$.

$2)$ to generalize the conditions for the random weights, staying always in frame of weighted form of \eqref{eq.KLP.1.8}. This way, we generalize both \cite[Th. 6.1]{konstantinides:passalidis:2025} and Theorem \ref{th.KLP.3.1} in this paper.

The rest of the paper is organized as follows. In section 2 we give some preliminary results for the distributions with heavy tails, with one or two dimensions and we introduce the class of regular variation in two dimensions. After depicting some dependence structures with some known results, a new dependence structure is introduced. Next, in section 3 we present the main result for the joint tail asymptotic behavior of the sums, together with some preliminary lemmas. In section 4 we study some closure properties of our two-dimensional classes, and we give some new results both in one and  two-dimensional cases. Finally in  section 5, similarly we provide the results for the randomly weighted sums as in case of dependence as also in case of GTAI structure, and we discuss the applications in ruin probabilities over finite time in a two-dimensional, discrete time, risk model with stochastic discount factors.

\section{Preliminaries} \label{sec.KP.2}

For two real numbers $x,\,y$ we denote $x^+:=\max\{x,\,0\}$, $x\wedge y:=\min\{x,\,y\}$ and $x \vee y := \max\{x,\,y\}$. With $ {\bf 1}_{A}$ we symbolized the indicator function on some event $A$. The vectors are denoted with bold script and ${\bf 0}$ denotes the origin of the axes. The joint distribution of two random variables $X,\,Y$ is given by ${\bf F}(x,\,y)=\PP[X\leq x\,,\;Y\leq y]$. As tail of the joint distribution we take the joint excess, namely ${\bf \bF}(x,\,y)={\bf \bF_1 }(x,\,y):= \PP[X> x, \, Y>y]$, where for sake of easiness, for any vector ${\bf b}=(b_1,\,b_2)$ we write ${\bf \bF_b}(x,\,y):= \PP[X>b_1\,x, \, Y>b_2\,y]$.

In what follows we use the asymptotic notations: for two uni-variate positive functions $f_1$ and $g_1$, with $f_1(x) \sim g_1(x)$, as $\xto$, we mean that
\beao
\lim_{\xto} \dfrac{f_1(x)}{g_1(x)}=1\,,
\eeao
with $f_1(x) \lesssim g_1(x)$, or $g_1(x) \gtrsim f_1(x)$, as $\xto$, we mean that
\beao
\limsup_{\xto} \dfrac{f_1(x)}{g_1(x)} \leq 1\,,
\eeao
with $f_1(x) =o[ g_1(x)]$, or $f_1(x) =o(1)\,g_1(x)$, as $\xto$, we mean that
\beao
\lim_{\xto} \dfrac{f_1(x)}{g_1(x)} =0\,,
\eeao
while with $f_1(x) =O[ g_1(x)]$, or $f_1(x) =O(1)\,g_1(x)$, as $\xto$, we mean that
\beao
\limsup_{\xto} \dfrac{f_1(x)}{g_1(x)} <\infty\,,
\eeao
and with $f_1(x) \asymp g_1(x)$, as $\xto$, we mean that
\beao
f_1(x)=O[g_1(x)]\, \qquad g_1(x)=O[f_1(x)]\,,
\eeao 
as $\xto$.

If there exist two positive two-variate functions $f_2$ and $g_2$, the corresponding limit relations come as extension of the one-variate with the restriction that the limits of two-variate functions hold for $x\wedge y \to \infty$. For example,  with $f_2(x,\,y) \sim g_2(x,\,y)$, as $x \wedge y \to \infty$, we mean that
\beao
\lim_{x \wedge y \to \infty} \dfrac{f_2(x,\,y)}{g_2(x,\,y)}=1\,,
\eeao
and $f_2(x,\,y) =o[g_2(x,\,y)]$, as $x \wedge y \to \infty$, if
\beao
\lim_{x \wedge y \to \infty} \dfrac{f_2(x,\,y)}{g_2(x,\,y)}=0\,.
\eeao

\subsection{Heavy-tailed Distributions}

The heavy-tailed distributions recently attract more interest of applied and theoretical probability community. For example in Risk Theory, Financial Mathematics, Risk Management, Branching Processes, L\'{e}vy Processes present the domains of applications with heavy tails, see for example  \cite{athreya:ney:1972}, \cite{asmussen:klueppelberg:1996}, \cite{asmussen:2003}, \cite{rolski:schmidli:schmidt:teugels:1999}, \cite{korshunov:2018}, \cite{foss:korshunov:palmowski:2024}. We give some of the most important classes of heavy-tailed distributions, with some of these classes to be used later, together with their properties.

We say that a distribution $F$ is heavy tailed, symbolically $F \in \mathcal{H}$ if for any $\vep >0$ holds
\beao
\int_{-\infty}^{\infty} e^{\vep \,y}\,F(dy) = \infty\,.
\eeao
We say that a distribution $F$ is long tailed, symbolically $F \in \mathcal{L}$ if for any (or equivalently, for some) $a>0$, holds
\beao
\lim_{\xto} \dfrac {\bF(x-a)}{\bF(x)} = 1\,.
\eeao 
It is well-known that, if $F \in \mathcal{L}$ then there exists some function $a\;:\;[0,\,\infty) \to [0,\,\infty)$ such that $a(x) \to \infty$, $a(x)=o(x)$ and $\bF[x\pm a(x)] \sim \bF(x)$, as $\xto$. Function  $a(x)$ is called insensitivity function for distribution $F$, see in \cite{konstantinides:2018} or in \cite{foss:korshunov:zachary:2013}. 

Let remind that for two random variables $X_1$ and $X_2$ with distributions $F_1$ and $F_2$ respectively, the distribution of the sum is defined as $F_{X_1+X_2}(x):=\PP[X_1+ X_2 \leq x]$ and its tail $\bF_{X_1+X_2}(x):=\PP[X_1 + X_2 >x]$.  If the random variables are independent then we denote $F_1*F_2$ instead of $F_{X_1+X_2}$, which is called convolution of $F_1$ and $F_2$. A distribution $F$ with support $\bbr_+$ is called subexponential, symbolically $F\in \mathcal{S}$, if
\beao
\lim_{\xto}\dfrac {\overline{F^{*n}}(x)}{\bF(x)} =n\,,
\eeao
for any (or equivalently, for some) integer $n\geq 2$. With $F^{*n}$ we depict the $n$-th convolution power of $F$ with itself. The classes $\mathcal{H}$, $\mathcal{L}$ and $\mathcal{S}$ were defined in \cite{chistyakov:1964}. 

If a distribution $F$ has support $\bbr$, then we say that it is subexponential on the real axis, if $F_+(x) \in \mathcal{S}$, where $F_+(x):=F(x)\,{\bf 1}_{\{x>0\}}$, see \cite{pakes:2004}. Therefore, if $X_1,\,\ldots,\,X_n$ are independent, identically distributed, real valued, random variables with common distribution $F \in \mathcal{S}$, we find relation \eqref{eq.KLP.1.1} and by elementary inequalities we get also  \eqref{eq.KLP.1.2}. This is a main reason to make class $\mathcal{S}$ very important, see for applications in risk theory in \cite{jiang:wang:chen:xu:2015}, \cite{li:2017}, \cite{yuan:lu:2023} among others.  

Class $\mathcal{D}$ was introduced in \cite{feller:1969} as extension of the regular variation. We say that a distribution $F$ belongs to the class of dominatedly varying distributions, symbolically $F \in \mathcal{D}$, if for any (or equivaletly, for some) $b \in (0,\,1)$ holds
\beao
\limsup_{\xto} \dfrac{\bF(b\,x)}{\bF(x)} < \infty\,.
\eeao
It is well-known that $\mathcal{D} \not\subseteq \mathcal{S}$, $\mathcal{S} \not\subseteq \mathcal{D}$ and $\mathcal{D}\cap \mathcal{S}\equiv \mathcal{D}\cap \mathcal{L}$, see \cite{goldie:1978}.

Another distribution class with heavy tails is class $\mathcal{C}$, of consistently varying distributions. We say that $F$ is consistently varying and we write $F \in \mathcal{C}$ if
\beao
\lim_{z \uparrow 1} \limsup_{\xto} \dfrac{\bF(z\,x)}{\bF(x)} =1\,.
\eeao

Further, we say that a distribution $F$ belongs to the class of regularly varying distributions, with index $\alpha>0$,  symbolically $F \in \mathcal{R}_{-\alpha}$, if holds
\beao
\lim_{\xto} \dfrac{\bF(t\,x)}{\bF(x)} =t^{-\alpha}\,,
\eeao 
for any $t>0$. It is well-known that $ \mathcal{R}:=\bigcup_{\alpha >0} \mathcal{R}_{-\alpha} \subsetneq\mathcal{C} \subsetneq \mathcal{D}\cap \mathcal{L} \subsetneq \mathcal{S} \subsetneq \mathcal{L}\subsetneq \mathcal{H}$, see for example \cite{konstantinides:2018}, \cite{leipus:siaulys:konstantinides:2023}.

 Next for any distribution with upper unbounded support, we bring up the upper and lower Matuszewska indexes $J_F^+$, $J_F^-$ respectively, introduced in \cite{matuszewska:1964}, which have important role on the characterization of heavy-tailed and related distributions
\beao
J_F^+ := -\lim_{v \rightarrow \infty }\dfrac{\log \bF_{*}(v)}{\log v }\,, \qquad J_F^- := -\lim_{v \rightarrow \infty }\dfrac{\log \bF^{*}(v)}{\log v }\,.
\eeao
where
\beao
\bF_{*}(v):=\liminf_{x\rightarrow \infty }\dfrac{\bF(v\,x)}{\bF(x)}\,, \qquad \bF^{*}(v):=\limsup_{x\rightarrow \infty }\dfrac{\bF(v\,x)}{\bF(x)}\,,
\eeao
for any $v>1$.
The following inequalities hold $0\leq J_F^- \leq J_F^+ \leq \infty$. We have the equivalence, $F\in \mathcal{D}$ if and only if $J_F^+ < \infty$ and, if $F\in \mathcal{R}_{-\alpha}$ for some $\alpha>0$, then $J_F^+=J_F^-=\alpha$.

\subsection{Two-dimensional Distribution Classes with Heavy Tails}

Several attempts appeared recently for the modeling of extreme events in multi-dimensional set up. The reason was the dependence structures among these events, that play crucial role on final outcome. Although the multivariate regular variation ($MRV$) is a well established multivariate extension of the regular variation to many dimensions, this is NOT the case for the rest distribution classes. The most popular attempts for multivariate distribution classes are referred to subexponential distributions and to long tailed distributions. In literature we find at least four different definitions to multivariate subexponential distributions, see \cite{cline:resnick:1992}, \cite{omey:2006}, \cite{samorodnitsky:sun:2016}, \cite{konstantinides:passalidis:2025}. From the first three, this by \cite{samorodnitsky:sun:2016}, seems as the strongest one. While the first three forms of the multivariate subexponentiality focus on multivariate linear single big jump principle, recall relation \eqref{eq.KLP.1.A.1}, the class in \cite{konstantinides:passalidis:2025} focus on multivariate non-linear single big jump principle,  namely on relation \eqref{eq.KLP.1.6}. 

The characterization of the bi-variate classes is based on the joint excess of random variables, as also on the marginal distributions. This makes the extension of one-dimensional classes more direct and in the same time permits arbitrary dependence structure among the components. Further, with the notation ${\bf a}=(a_1,\,a_2)> (0,\,0)$, we have in mind that  ${\bf a} \in [0,\,\infty)^2 \setminus \{{\bf 0}\}$, except it is referred differently. The following classes were introduced in \cite{konstantinides:passalidis:2025}.

Let consider initially the two-dimensional class of long tailed distributions.
We say that a random pair $(X,\,Y)$ with marginal distributions $F$, $G$, belongs to the class of bi-variate long tailed, symbolically ${\bf F}:=(F,\,G) \in \mathcal{L}^{(2)}$, if 
\begin{enumerate}

\item

$F \in \mathcal{L}$, $G \in \mathcal{L}$,

\item

\beao
\lim_{x\wedge y \to \infty} \dfrac{{\bf \bF}_1(x-a_1,\,y-a_2)}{{\bf \bF}_1(x,\,y)} = \lim_{x\wedge y \to \infty} \dfrac{\PP[X>x-a_1\,,\;Y>y-a_2]}{\PP[X>x\,,\;Y>y]} =1\,,
\eeao
for some, (or equivalently, for any) ${\bf a}=(a_1,\,a_2)> (0,\,0)$.
\end{enumerate} 
 
As in one-dimensional case there exists insensitivity function, in two-dimensional case also exists $a$-jointly insensitivity function. Concretely, for a pair of two non-decreasing functions $a_F(x),\,a_G(y) >0$, for any $x,\,y>0$, where $a_F(x) \to \infty$ and $a_F(x)=o(x)$ and $a_G(y) \to \infty$ and $a_G(y)=o(y)$, as $\xto$ and $y \to \infty$ respectively, the  two-dimensional distribution ${\bf F} = (F,\,G)$, with infinite right endpoints for both marginal distributions, is called $(a_F,\,a_G)$-joint insensitivity distribution, if the following is true
\beam \label{eq.KLP.2.3} \notag
&&\sup_{|a_1|\leq a_F(x),  |a_2|\leq a_G(y)} \left|\PP[X>x-a_1, Y>y-a_2]-\PP[X>x, Y>y]\right| \\[2mm] 
&&\qquad \qquad =o\left( \PP[X>x, Y>y]\right)\,,
\eeam
as $x\wedge y \to \infty$. From \cite[Lem. 2.1]{konstantinides:passalidis:2025} we find that if 
$(F,\,G) \in \mathcal{L}^{(2)}$, then there exists some bi-variate function $(a_F(x),\,a_G(y))$, such that relation \eqref{eq.KLP.2.3} holds. Furthermore, this means that
\beao
\lim_{x\wedge y \to \infty} \dfrac{\PP[X>x\pm a_F(x)\,,\;Y>y\pm a_G(y)]}{\PP[X>x\,,\;Y>y]} =1\,,
\eeao
but this function is NOT making $a$-insensitive the marginal distributions.

We say that the random pair $(X,\,Y)$, with marginal distributions $F$ and $G$ respectively, belongs to the class of bi-variate subexponential distributions, symbolically $(F,\,G) \in \mathcal{S}^{(2)}$, if
\begin{enumerate}
\item
$F \in \mathcal{S}$ and $G \in \mathcal{S}$.
\item
$(F,\,G) \in \mathcal{L}^{(2)} $.
\item
It holds
\beam \label{eq.KP.13}
\lim_{x\wedge y \to \infty} \dfrac{\PP[X_1+X_2>x\,,\;Y_1+Y_2>y]}{\PP[X>x\,,\;Y>y]}=2^2\,,
\eeam
where $(X_1\,,\;Y_1)$ and $(X_2\,,\;Y_2)$ are independent and identically distributed  copies of $(X\,,\;Y)$. 
\end{enumerate}

\bre
In case of $d$-variate distribution relation \eqref{eq.KP.13} becomes 
\beam \label{eq.KP.13b}
\lim_{x_1\wedge \cdots \wedge x_d \to \infty} \dfrac{\PP[X_{1,1}+X_{1,2}>x_1\,,\;\ldots,\,X_{d,1}+X_{d,2}>x_d]}{\PP[X_{1,1}>x_1\,,\;\ldots,\,X_{d,1}>x_d]}=2^d\,.
\eeam
From \eqref{eq.KP.13b} we can easily see that the subexponentiality by \cite{konstantinides:passalidis:2025} satisfies the multivariate non-linear single big jump principle in \eqref{eq.KLP.1.6}, in what fails also the multidimensional subexponentiality from \cite{samorodnitsky:sun:2016}, but also the well-established class of standard $MRV$, since these classes are constructed in way that they satisfy the multivariate linear single big jump by relation \eqref{eq.KLP.1.A.1}.
\ere

Another bi-variate distribution class is the $\mathcal{D}^{(2)}$ of the bi-variate dominatedly varying distributions. For a random pair $(X,\,Y)$ with marginal distributions $F$, $G$, we say that belongs to the class of bi-variate dominatedly varying distributions and we write $(F,\,G) \in \mathcal{D}^{(2)}$, if 
\begin{enumerate}
\item
$F \in \mathcal{D}$, $G \in \mathcal{D}$,
\item
\beao
\limsup_{x\wedge y \to \infty} \dfrac{{\bf \bF_b}(x,\,y)}{{\bf \bF_1}(x,\,y)} = \limsup_{x\wedge y \to \infty} \dfrac{\PP[X>b_1\,x\,,\;Y>b_2\,y]}{\PP[X>x\,,\;Y>y]} < \infty\,,
\eeao
for any, (or equivalently, for some) ${\bf b}=(b_1,\,b_2) \in (0,\,1)^2$.
\end{enumerate} 

Let us denote $(\mathcal{D}\cap \mathcal{L})^{(2)}:=\mathcal{D}^{(2)} \cap \mathcal{L}^{(2)}$. 

Next, we say that the random pair $(X,\,Y)$ with marginal distributions $F$, $G$, belongs to the class of bi-variate consistently varying distributions, symbolically $(F,\,G) \in \mathcal{C}^{(2)}$, if 
\begin{enumerate}

\item
$F \in \mathcal{C}$, $G \in \mathcal{C}$,

\item
\beam \label{eq.KLP.2.5}
\lim_{{\bf z} \uparrow {\bf 1}}\limsup_{x\wedge y \to \infty} \dfrac{{\bf \bF_z}(x,\,y)}{{\bf \bF_1}(x,\,y)} =\lim_{{\bf z} \uparrow {\bf 1}} \limsup_{x\wedge y \to \infty} \dfrac{\PP[X>z_1\,x\,,\;Y>z_2\,y]}{\PP[X>x\,,\;Y>y]} =1\,,
\eeam
with ${\bf z}=(z_1,\,z_2)$ and ${\bf 1}=(1,\,1)$.
\end{enumerate} 
It was proved that $\mathcal{C}^{(2)}\subsetneq (\mathcal{D}\cap \mathcal{L})^{(2)}$, see \cite[Th. 2.1]{konstantinides:passalidis:2025}.

Now we introduce a new distribution class, which represents a new kind of two-variate regular variation. Except the definition we comment about the advantages and disadvantages in relation with the classical multivariate regular variation, introduced in \cite{dehaan:resnick:1982}.

\bde \label{def.KLP.2.1}
We say that the random pair $(X,\,Y)$ with marginal distributions $F$, $G$ respectively, follows two-dimensional regular variation, symbolically ${\bf F} \in \mathcal{R}_{(-\alpha_1,\,-\alpha_2)}^{(2)}$, with $0<\alpha_1, \alpha_2 < \infty$, if 
\begin{enumerate}
\item
$F \in \mathcal{R}_{-\alpha_1}$, $G \in \mathcal{R}_{-\alpha_2}$,
\item
\beam \label{eq.KLP.2.6}
\lim_{x\wedge y \to \infty} \dfrac{{\bf \bF_t}(x,\,y)}{{\bf \bF_1}(x,\,y)} = \lim_{x\wedge y \to \infty} \dfrac{\PP[X>t_1\,x\,,\;Y>t_2\,y]}{\PP[X>x\,,\;Y>y]} = t_1^{-\alpha_1}\,t_2^{-\alpha_2}\,,
\eeam
for any ${\bf t}=(t_1,\,t_2)$, with $t_1,\,t_2> 0$.
\end{enumerate} 
\ede

In this definition, for the case $\alpha_1=\alpha_2$, we say that we have a typical, two-dimensional, regular variation, while in opposite case  $\alpha_1 \neq \alpha_2$, we have a non-typical,  two-dimensional, regular variation.

\bre \label{rem.KLP.2.1}
Let restrict ourselves to two-dimensional case, say the $BRV$, for the depiction of the multivariate regular variation. If the normalizing functions are
\beao
U_X(x)=\left(\dfrac 1{\bF} \right)^{\leftarrow} (x)\,, \qquad U_Y(x)=\left(\dfrac 1{\bG} \right)^{\leftarrow} (x)\,,
\eeao
where $f^{\leftarrow}$ represents the c\'{a}gl\'{a}d inverse of function $f$, for $x>0$, then the bi-variate regular variation, symbolically $(X\,,\;Y) \in BRV_{-\alpha_1,\,-\alpha_2}(F,\,G,\,\nu)$, is defined as
\beao
\lim_{\xto} x\,\PP\left[\left(\dfrac X{U_X(x)}\,,\;\dfrac Y{U_Y(x)}\right)\in \bbb \right]= \nu(\bbb)\,,
\eeao
for any Borel set $\bbb \subsetneq [0,\,\infty]^2 \setminus \{ {\bf 0} \}$, and the limit measure $\nu$ is non-degenerate. Furthermore, it follows that the distributions $F$ and $G$ are regularly varying and the measure $\nu$ is homogeneous in the sense of $\nu(\bbb^{\lambda})=\lambda^{-1}\,\nu(\bbb)$, where $\bbb^{\lambda} = \left\{\left(\lambda^{1/\alpha_1}\,\kappa\,,\;\lambda^{1/\alpha_2}\,\mu\right)\;:\;(\kappa,\,\mu) \in \bbb \right\}$, see in \cite[Lem. 5.1]{tang:yang:2019}, where $\alpha_1,\,\alpha_2$ are the regular variation indexes of $F,\,G$ respectively. If $U_X:=U_Y$, we have standard $BRV$.
\ere

The class $MRV$ attracted the interest of several branches of applied probability and statistics, see for example \cite{basrak:davis:mikosch:2002}, \cite{buraczewski:damek:mikosch:2016}, \cite{resnick:2007}.

Recently, risk theory gained popularity, as for example in \cite{chen:yang:2019}, \cite{tang:yang:2019} in non-standard $BRV$, and for example in \cite{cheng:konstantinides:wang:2024}, \cite{konstantinides:li:2016}, \cite{yang:su:2023}, in standard $MRV$. 

Although, class $BRV$ is well-established in standard case, there are several difficulties, coming from the non-standard case, that appear because of the variety of the normalizing functions. Furthermore, the main issue is that it can NOT handle the multivariate non-linear single big jump by relation \eqref{eq.KLP.1.6}, which represents the main question we examine in this paper. However, we should stress that since $\mathcal{R}_{(-\alpha_1,\,-\alpha_2)}^{(2)} \subsetneq \mathcal{S}^{(2)}$, from relations \eqref{eq.KP.13} and \eqref{eq.KLP.1.A.1}, follows directly that
\beao
\mathcal{R}_{(-\alpha_1,\,-\alpha_2)}^{(2)} \bigcap BRV_{-\alpha_1,\,-\alpha_2}(F,\,G,\,\nu) = \emptyset\,, 
\eeao
for any $0<\alpha_1, \alpha_2 < \infty$. The class $\mathcal{R}_{(-\alpha_1,\,-\alpha_2)}^{(2)}$ tries to surpass this obstacle, but with loss of some dependence cases, as the asymptotic dependence, see Proposition \ref{prop.KP.2.1}, below in this section. 

Let us provide some examples for the building of the class $\mathcal{R}_{(-\alpha_1,\,-\alpha_2)}^{(2)}$. At first, it is obvious that if $F\in \mathcal{R}_{-\alpha_1}$ and $G\in \mathcal{R}_{-\alpha_2}$, with $a_1,\,a_2 >0$, with their distributions stemming from independent random variables, then it  holds ${\bf F} \in \mathcal{R}_{(-\alpha_1,\,-\alpha_2)}^{(2)}$. An interesting case of dependence presents the strongly asymptotic independence, symbolically $SAI$, see \cite{li:2018a}, \cite{li:2018b}. Without loss of generality, for sake of simplicity with less conditions, we consider two non-negative random variables $X,\,Y$ with distributions $F,\,G$ respectively. We say that the $X,\,Y$ satisfy the $SAI$ condition, if there exists a constant $C > 0$ such that
\beao
\PP[X>x\,,\;Y>y]=[C+o(1)]\,\bF(x)\,\bG(y)\,,
\eeao
as $x\wedge y \to \infty$. 

Under the previous conditions, if $F \in \mathcal{R}_{-\alpha_1}$ and $G \in \mathcal{R}_{-\alpha_2}$, then for any $t_1,\, t_2 >0$ it holds
\beao
\lim_{x\wedge y \to \infty} \dfrac{{\bF_t}(x,\,y)}{{\bF_1}(x,\,y)}=\dfrac{C\,\bF(t_1\,x)\,\bG(t_2\,y)}{C\,\bF(x)\,\bG(y)} = t_1^{-\alpha_1}\,t_2^{-\alpha_2}\,,
\eeao
from where, we obtain ${\bf F} \in \mathcal{R}_{(-\alpha_1,\,-\alpha_2)}^{(2)}$.

Next, we present some relations with respect to the characterization of closure properties in a two-dimensional class, as also of joint max-sum equivalence, as can be found in \cite{konstantinides:passalidis:2025}. 
\begin{enumerate}
\item
Closure property with respect to sum. Let $X_1,\,X_2,\,Y_1,\,Y_2$ random variables with distributions $F_1,\,F_2,\,G_1,\,G_2$ respectively. If the condition $F_1,\,F_2,\,G_1,\,G_2 \in \mathcal{B}$ holds, for any $k,\,l \in \{1,\,2\}$ holds $(F_k,\,G_l) \in \mathcal{B}^{(2)}$, and also holds $(F_{X_1+X_2},\,G_{Y_1+Y_2}) \in \mathcal{B}^{(2)}$, then we say that $\mathcal{B}^{(2)}$ is closed with respect to the sum, where $\mathcal{B}^{(2)}$ represents some two-dimensional version of some class $\mathcal{B}$.
\item
Closure property with respect to convolution product. If $X$, $Y$ are random variables with distributions $F$, $G$ respectively, it holds $(F,\,G)\in \mathcal{B}^{(2)}$ and $(\Theta,\,\Delta)$ is a random pair, independent of $(X,\,Y)$, then we say that the two-dimensional class $\mathcal{B}^{(2)}$ is closed with respect to convolution product if $(\Theta\,X,\,\Delta\,Y) \in \mathcal{B}^{(2)}$.
\item
Joint max-sum equivalence.  Let $X_1,\,X_2,\,Y_1,\,Y_2$ be random variables, then we say that they are joint max-sum equivalent if it holds $\PP[ X_1+X_2 >x\,,\;Y_1+Y_2 >x] \sim \sum_{k=1}^2 \sum_{l=1}^2 \PP[X_k > x\,,\;Y_k >y ]$, as $x\wedge y \to \infty$. 
\end{enumerate}

\bre \label{rem.KLP.2.2}
In one dimensional set up the closure properties are well studied with respect to convolution, convolution product, minimums and maximums. For some papers on this topic see  \cite{cline:samorodnitsky:1994}, \cite{tang:2006}, \cite{cui:wang:2020}, \cite{leipus:siaulys:2020}. Furthermore, a detailed account of closure properties for heavy tailed distributions we refer to \cite{leipus:siaulys:konstantinides:2023}. In multidimensional set up, there are only a few attempts, mostly with respect to convolution product of $MRV$, see \cite{basrak:davis:mikosch:2002b} and \cite{fougeres:mercadier:2012}. An introduction of random vectors with heavy tails in relation with closure properties can be found in \cite{konstantinides:passalidis:2024b}, \cite{das:fasenhartmann:2023} and \cite{konstantinides:passalidis:2024g}.
\ere

\subsection{Dependence Modeling}

Now, we introduce the dependence structures, which are useful in our results. We remind two dependencies, that apply on only single sequence of random variables and next we provide another two structures, which apply on double sequences of random variables.

For a sequence of real valued random variables $X_1,\,\ldots,\,X_n$ with distributions $F_1,\,\ldots,\,F_n$ respectively, we say that they are pairwise quasi-asymptotically independent, symbolically $QAI$, if for any pair $i,\,j =1,\,\ldots,\,n$, with $i \neq j$ it holds 
\beam \label{eq.KLP.2.9}
\lim_{\xto}\dfrac{\PP[ |X_i| >x\,,\;X_j >x]}{\bF_i(x) + \bF_j(x) } = 0\,.
\eeam
For the same sequence of random variables, we say that they are Tail Asymptotic Independent, symbolically $TAI$, and in some works named as strong quasi-asymptotically independent, if for any  pair $i,\,j =1,\,\ldots,\,n$, with $i \neq j$ holds the limit
\beam \label{eq.KLP.2.10}
\lim_{x_i \wedge x_j \to \infty}\PP[ |X_i| >x_i\;|\; X_j >x_j] = 0\,.
\eeam
From relations \eqref{eq.KLP.2.9} and \eqref{eq.KLP.2.10} we see easily the inclusion $TAI \subsetneq QAI$. The dependencies $QAI$ and $TAI$ were introduced in \cite{chen:yuen:2009} and \cite{geluk:tang:2009} respectively, where the asymptotic formula \eqref{eq.KLP.1.1} was proved for these two dependencies in the distribution classes $\mathcal{C}$ and $\mathcal{D}\cap \mathcal{L}$ respectively. These dependence structures are included  in the concept of asymptotic independence, see for example \cite{maulik:resnick:2004}, and used for randomly weighted sums, or for generalized moments of randomly weighted sums, see \cite{wang:2011}, \cite{li:2013}, \cite{dirma:nakiliuda:siaulys:2023}, \cite{konstantinides:leipus:passalidis:siaulys:2025}.

The following proposition shows that class $\mathcal{R}_{(-\alpha_1,\,-\alpha_2)}^{(2)}$ is restricted in the $QAI$ structure, for the random pair $(X,\,Y)$.

\bpr \label{prop.KP.2.1}
If $(X,\,Y) \in \mathcal{R}_{(-\alpha_1,\,-\alpha_2)}^{(2)}$, then $X$ and $Y$ are $QAI$.
\epr

\pr~
We use the method of contradiction. Let assume that $X$ and $Y$ are not quasi-asymptotically independent. Then, at least one of the following two relations does not hold
\beam \label{eq.KP.2.*}
\lim_{\xto} \dfrac{\PP[X>x,\,|Y|>x]}{\bF(x)+\bG(x)} =0\,,\qquad \lim_{\xto} \dfrac{\PP[|X|>x,\,Y>x]}{\bF(x)+\bG(x)} =0\,.
\eeam
Without loss of generality, let $\alpha_1 \leq \alpha_2$ and let assume that the first relation in \eqref{eq.KP.2.*} does not hold. Note also that
\beao
\lim_{\xto} \dfrac{\PP[X>x,\,|Y|>x]}{\bF(x)+\bG(x)} \leq 1\,.
\eeao
Hence, we find
\beao
\limsup_{\xto} \dfrac{\PP[X>x,\,|Y|>x]}{\bF(x)+\bG(x)} =:c \in (0,\,1]\,.
\eeao
Thus, there exists a sequence $\{x_n,\;n\in \bbn\}$, such that the convergence $\lim_{\nto} x_n= \infty$ and 
\beao
\lim_{\nto} \dfrac{\PP[X>x_n,\,|Y|>x_n]}{\bF(x_n)+\bG(x_n)}=c\,,
\eeao
hold. We choose now $t$ small enough, such that $1/c < t^{-\alpha_1}$, which implies $t \in (0,\,1)$. Then, we obtain
\beao
&&\lim_{\nto} \dfrac{\PP[ X >t\,x_n,\; Y >t\,x_n]}{\PP[ X >x_n,\;Y>x_n]}\leq \lim_{\nto} \dfrac{\PP[ X >t\,x_n,\; |Y|>t\,x_n]}{\PP[ X >x_n,\;Y>x_n]}\\[2mm]
&&=\lim_{\nto} \dfrac{\PP[ X >t\,x_n,\; |Y| >t\,x_n]}{\bF(t\,x_n)+\bG(t\,x_n)}\, \dfrac{\bF(t\,x_n)+\bG(t\,x_n)}{\bF(x_n)+\bG(x_n)}\, \dfrac{\bF(x_n)+\bG(x_n)}{\PP[ X >x_n,\;Y>x_n]}\\[2mm]
&&\leq \lim_{\nto} \max\left\{ \dfrac{\bF(t\,x_n)}{\bF(x_n)}\,,\;\dfrac{\bG(t\,x_n)}{\bG(x_n)} \right\} \dfrac{\bF(x_n)+\bG(x_n)}{\PP[ X >x_n,\;Y>x_n]}=\dfrac{t^{-\alpha_2}}c < t^{-\alpha_1}\,t^{-\alpha_2}\,.
\eeao 
This implies that
\beao
&&\lim_{\nto} \dfrac{\PP[ X >t\,x,\; Y >t\,x]}{\PP[ X >x,\;Y>x]} \neq t^{-\alpha_1}\,t^{-\alpha_2}\,,
\eeao
which contradicts with $(X,\,Y) \in \mathcal{R}_{(-\alpha_1,\,-\alpha_2)}^{(2)}$.
~\halmos

\bre \label{rem.KP.2.4*}
Another observation to class $\mathcal{R}_{(-\alpha_1,\,-\alpha_2)}^{(2)}$ is the following. Let us write
\beao
f(x,\,y):= \dfrac {\PP[X>x,\,Y>y]}{\PP[X>x]\,\PP[Y>y]}\,,
\eeao
then, if $(X,\,Y) \in \mathcal{R}_{(-\alpha_1,\,-\alpha_2)}^{(2)}$, with $\alpha_1,\,\alpha_2 >0$, the function $f(x,\,y)$ has some kind of 'bi-variate slow variation', in the sense that for any $t_1,\,t_2>0$ it holds
\beao
\lim_{x \wedge y \to \infty} \dfrac{f(t_1\,x,\,t_2\,y)}{f(x,\,y)}=\lim_{x \wedge y \to \infty} \dfrac{\dfrac {\PP[X>t_1\,x,\,Y>t_2\,y]}{\PP[X>t_1\,x]\,\PP[Y>t_2\,y]}}{\dfrac {\PP[X>x,\,Y>y]}{\PP[X>x]\,\PP[Y>y]}}=1\,,
\eeao
hence we find
\beao
\dfrac {\PP[X>t_1\,x,\,Y>t_2\,y]}{\PP[X>x\,,\;Y>y]}=\dfrac{f(t_1\,x,\,t_2\,y)\,\bF(t_1\,x)\,\bG(t_2\,y)}{f(x,\,y)\,\bF(x)\,\bG(y)} \to t^{-\alpha_1}\,t^{-\alpha_2}\,,
\eeao
as $x \wedge y \to \infty$. By this is implied a generalization of the $SAI$ dependence structure. Indeed, we get $\PP[X>x,\,Y>y]=f(x,\,y)\,\bF(x)\,\bG(y)$, to become $SAI$, when the function $f(x,\,y)$ takes a constant value, and the last relation is not necessarily an asymptotic one.
\ere

In the dependence structures with two sequences of random variables, we want to model simultaneously the dependence among the terms of each sequences and the dependence between the two sequences. Let two sequences of real valued, random variables  $\{X_1,\,\ldots,\,X_n\}$  and $\{Y_1,\,\ldots,\,Y_m\}$. We say that they are Generalized Tail Asymptotic Independent, symbolically $GTAI$, if 
\beao
\lim_{x_i \wedge x_k \wedge y_j \to \infty}\PP[ |X_i| >x_i\;|\; X_k >x_k,\;Y_j>y_j] = 0\,,
\eeao  
for any $i,\,k =1,\,\ldots,\,n$, $j=1,\ldots, \,m$, with $i \neq k$, and if
\beao
\lim_{x_i \wedge y_k \wedge y_j \to \infty}\PP[ |Y_j| >y_j\;|\;Y_k>y_k,\; X_i >x_i] = 0\,,
\eeao  
for any $j,\,k =1,\,\ldots,\,m$, $i=1,\ldots, \,n$, with $j \neq k$.

This kind of dependence structure indicates that the probability of happening three or more extreme events is negligible in comparison with the probability of happening two extreme events, one in each sequence. Furthermore, we see that if the two sequences $\{X_1,\,\ldots,\,X_n\}$  and $\{Y_1,\,\ldots,\,Y_m\}$ are independent, then each one of them has $TAI$ dependence structure. From the other hand side, under $GTAI$ structure, if the $\{(X_i,\,Y_i)\,,\;i \in \bbn\}$ are independent random pairs, then each pair $(X_i,\,Y_i)$ has $TAI$ components. The $GTAI$ dependence structure can be found in \cite{konstantinides:passalidis:2023}.

Next, we introduce a new dependence structure between the two sequences, that generalizes the $GTAI$, but remains in the same spirit.

\bde \label{def.KLP.2.2}
Let two sequences of real valued, random variables  $\{X_1,\,\ldots,\,X_n\}$  and $\{Y_1,\,\ldots,\,Y_m\}$. We say that they are Generalized Quasi Asymptotically Independent, symbolically $GQAI$, if
\beao
\lim_{x \wedge y \to \infty} \dfrac{\PP[ |X_i| >x\ , \;X_k >x,\;Y_j>y]}{\PP[ X_i >x,\;Y_j>y] +\PP[ X_k >x,\;Y_j>y]}= 0\,,
\eeao 
for any $i,\,k =1,\,\ldots,\,n$, $j=1,\,\ldots,\,m$ with $i \neq k$, and if
\beao
\lim_{x \wedge y \to \infty}\dfrac{\PP[ |Y_j| >y, \; Y_k >y,\;X_i>x]}{\PP[X_i>x,\; Y_j >y] +\PP[ X_i >x,\;Y_k>y]}= 0\,,
\eeao 
for any $j,\,k =1,\,\ldots,\,m$, $i=1,\,\ldots,\,n$ with $j \neq k$.
\ede 

\bre \label{rem.KLP.2.3}
We can observe that $GTAI \subsetneq GQAI$, and if the two sequences are independent, then for any $i,\,k =1,\,\ldots,\,n$, $j=1,\,\ldots,\,m$ with $k\neq i$ we obtain
\beao
0&\sim& \dfrac{\PP[ |X_i| >x,\; X_k >x,\;Y_j>y]}{\PP[ X_i >x,\;Y_j>y] +\PP[ X_k >x,\;Y_j>y]}\\[2mm]
&=& \dfrac{\PP[ |X_i| >x,\; X_k >k]\,\PP[Y_j >y]}{(\PP[X_i>x] +\PP[ X_i >x])\,\PP[Y_j >y]}=\dfrac{\PP[ |X_i| >x,\; X_k >k]}{\bF_i(x) + \bF_k(x)} \sim 0\,,
\eeao
as $x \wedge y\to \infty$, hence, the $\{X_1,\,\ldots,\,X_n\}$ have the  $QAI$ dependence structure and similarly the $\{Y_1,\,\ldots,\,Y_m\}$ have also the  $QAI$ dependence structure.

Further, we see that the dependencies $GTAI$ and $ GQAI$ contain the complete independence as a special case, namely when the two sequences are independent and both sequences have independent terms.
\ere

The main target in this paper is the estimation of the asymptotic expressions \eqref{eq.KLP.1.6} and \eqref{eq.KLP.1.7}, in first case of  $GQAI$ dependence in distribution class $\mathcal{C}^{(2)}$, and in second case of $GTAI$ dependence in distribution class $(\mathcal{D} \cap \mathcal{L})^{(2)}$, for not necessarily common multitude of terms, namely for $n \neq m$.

\section{Joint Tail Behavior of Random Sums} \label{sec.KP.3}

In this section we show that the insensitivity property of the joint tail distribution of the sums, with respect to $GQAI$ dependence in the frame of distribution class $\mathcal{C}^{(2)}$. 

\bth \label{th.KLP.3.1}
Let two sequences of real valued, random variables  $\{X_1,\,\ldots,\,X_n\}$, with distributions $F_1,\,\ldots,\,F_n$  and $\{Y_1,\,\ldots,\,Y_m\}$, with distributions $G_1,\,\ldots,\,G_m$ respectively. If both sequences $X_1,\,\ldots,\,X_n,\,Y_1,\,\ldots,\,Y_m  \in GQAI$, and $(X_i,\,Y_j) \in \mathcal{C}^{(2)}$, for any $i=1,\,\ldots,\,n$ and $j=1,\,\ldots,\,m$, then
\beam \label{eq.KLP.3.1} \notag
\PP[ S_n >x,\,T_m>y]&\sim& \PP\left[ \bigvee_{i=1}^n S_i >x,\,\bigvee_{j=1}^m T_j>y\right] \sim \sum_{i=1}^n \sum_{j=1}^m \PP[ X_i >x,\,Y_j>y]\\[2mm]
&=:&S(x,\,y)\,,
\eeam  
as $x \wedge y\to \infty$.
\ethe 

\pr~
We start with the upper inequality
\beam \label{eq.KLP.3.2}
\PP [ S_n >x,\;T_m>y] \leq [1+o(1)] S(x,\,y)\,,
\eeam  
as $x \wedge y\to \infty$. Let choose arbitrarily a constant $\vep \in (0,\,1)$, then we obtain
\beam \label{eq.KLP.3.3}
&&\PP[ S_n >x,\;T_m>y]\leq  \PP\left[ \bigcup_{i=1}^n \left\{X_i >(1-\vep)\,x\right\},\;T_m >y \right]\\[2mm] \notag
&& \qquad \qquad  + \PP\left[S_n >x,\; T_m>y,\;\bigcap_{i=1}^n \left\{X_i \leq (1-\vep)\,x\right\} \right]=: I_1(x,\,y)+I_2(x,\,y) \,.
\eeam  
For the first term we find
\beam \label{eq.KLP.3.4}
&&I_1(x,\,y)\leq  \PP\left[ \bigcup_{i=1}^n \left\{X_i >(1-\vep)\,x\right\},\;
\bigcup_{j=1}^m \left\{Y_j >(1-\vep)\,y\right\} \right]\\[2mm] \notag
&& + \PP\left[ \bigcup_{i=1}^n \left\{X_i >(1-\vep)\,x\right\},\; T_m>y,\;\bigcap_{j=1}^m \left\{Y_j \leq (1-\vep)\,y\right\} \right]=:I_{11}(x,\,y) +I_{12}(x,\,y) \,.
\eeam  
Next, we estimate $I_{11}(x,\,y)$
\beam \label{eq.KLP.3.5}
I_{11}(x,\,y) \leq \sum_{i=1}^n \sum_{j=1}^m \PP\left[ X_i >(1-\vep)\,x,\;
Y_j >(1-\vep)\,y \right] \,.
\eeam 
For term $I_{12}(x,\,y) $ we obtain
\beam \label{eq.KLP.3.6} \notag
&& I_{12}(x,\,y) =  \PP\left[ \bigcup_{i=1}^n \left\{X_i >(1-\vep)\,x\right\},\; T_m>y,\; \bigcup_{l=1}^m \left\{Y_l >\dfrac ym \right\},\; \bigcap_{j=1}^m \left\{Y_j \leq (1-\vep)\,y\right\} \right]\\[2mm] 
&& \leq \sum_{i=1}^n \sum_{l=1}^m  \PP\left[ X_i >(1-\vep)\,x,\;
Y_l >\dfrac ym,\;T_m-Y_l > \vep\,y \right] \\[2mm] \notag
&&\leq \sum_{i=1}^n \sum_{l\neq j=1}^m \PP\left[ X_i >(1-\vep)\,x,\;
Y_l >v(y,\,m),\;Y_j > v(y,\,m) \right] =  \\[2mm] \notag
&&o\Bigg( \sum_{i=1}^n \sum_{l\neq j=1}^m\PP\bigg[ X_i >(1-\vep) x,\; Y_j > v(y,\,m) \bigg]  + \PP\left[ X_i >(1-\vep)\, x,\; Y_l > v(y,\,m) \right]\Bigg) \\[2mm] \notag
&&=o[S(x,\,y)]\,,
\eeam
where $v(y,\,m):=(y/m) \wedge [({\vep\,y})/({m-1})]$, as $x \wedge y\to \infty$, where in the last step we used that 
\beam \label{eq.K.3.7}
(X_i,\,Y_j) \in \mathcal{C}^{(2)} \subsetneq \mathcal{D}^{(2)}\,,
\eeam 
and in the pre-last step we take into account that the $\{X_1,\,\ldots,\,X_n\}$ and $\{Y_1,\,\ldots,\,Y_m\}$ belong to $GQAI$. Therefore, from \eqref{eq.KLP.3.5} and \eqref{eq.KLP.3.6} and by the arbitrary choice of $\vep > 0$, taking into consideration the relation \eqref{eq.KLP.3.4}, we finally obtain the following upper asymptotic bound of the $I_{1}(x,\,y)$
\beam \label{eq.KLP.3.7}
I_{1}(x,\,y) \leq[1+o(1)] \sum_{i=1}^n \sum_{j=1}^m \PP\left[ X_i >x,\;
Y_j >y \right] \,,
\eeam 
as $x \wedge y\to \infty$.

Now we consider the estimation of term $I_{2}(x,\,y)$
\beam \label{eq.KLP.3.8}
&&I_{2}(x,\,y)  \leq  \PP\left[ S_n>x,\;\bigcap_{i=1}^n \left\{X_i \leq(1-\vep)\,x\right\},\; \bigcup_{j=1}^m \left\{Y_j > (1- \vep)\,y \right\}\right]\\[2mm]  \notag
&&\qquad \qquad +\PP\left[  S_n>x,\;T_m > y,\;\bigcap_{i=1}^n \left\{X_i \leq(1-\vep)\,x\right\},\;\bigcap_{j=1}^m \left\{Y_j \leq (1-\vep)\,y\right\} \right]\\[2mm]  \notag
&&=: I_{21}(x,\,y) +I_{22}(x,\,y).
\eeam
But for the first term $I_{21}(x,\,y)$ holds
\beao
&&I_{21}(x,\,y)=\PP\left[ S_n>x,\;\bigcap_{i=1}^n \left\{X_i \leq(1-\vep)\,x\right\},\; \bigcup_{j=1}^m \left\{Y_j > (1- \vep)\,y \right\},\;\bigcup_{k=1}^n \left\{X_k >\dfrac xn \right\}\right]\\[2mm] 
&&\leq \sum_{k=1}^n \sum_{j=1}^m  \PP\left[ X_k >\dfrac xn,\;
S_n - X_k >\vep\,x ,\;Y_j > (1-\vep)\,y \right]\\[2mm] 
&&\leq \sum_{i\neq k=1}^n \sum_{j=1}^m  \PP\left[ X_i >v(x,\,n),\;
 X_k >v(x,\,n),\;Y_j > (1-\vep)\,y \right]= \\[2mm] 
&& o\Bigg[  \sum_{i=1}^n \sum_{j=1}^m\left( \PP\left[ X_i >v(x,\,n),\; Y_j > (1-\vep)\,y \right] + \PP\left[ X_k >v(x,\,n),\; Y_j > (1-\vep)\,y \right]\right)\Bigg]\,,
\eeao
as $x \wedge y\to \infty$. Hence,
\beam \label{eq.KLP.3.9}
I_{21}(x,\,y) = o\left[S(x,\,y)\right]\,,
\eeam
as $x \wedge y\to \infty$, where we used again \eqref{eq.K.3.7}
and that $X_1,\,\ldots,\,X_n,\,Y_1,\,\ldots,\,Y_m \in GQAI$.

For the $I_{22}(x,\,y)$ we find
\beao
&&I_{22}(x,\,y)=\PP\Bigg[ S_n>x,\;T_m>y,\;\bigcap_{i=1}^n \left\{X_i \leq(1-\vep)\,x\right\},\;\\[2mm]
&&\qquad \qquad  \bigcap_{j=1}^m \left\{Y_j > (1- \vep)\,y \right\}, \bigcup_{k=1}^n \left\{X_k >\dfrac xn \right\},\;\bigcup_{l=1}^m \left\{Y_l > \dfrac ym \right\} \Bigg]\\[2mm] 
&&\leq \sum_{k=1}^n \sum_{l=1}^m  \PP\left[ X_k >\dfrac xn,\;
S_n - X_k >\vep\,x ,\;Y_l > \dfrac ym,\;T_m-Y_l> \vep\,y \right] \\[2mm] 
&&\leq\sum_{i\neq k=1}^n \sum_{j\neq l=1}^m \PP\left[ X_i >v(x,\,n),\;X_k >v(x,\,n),\;Y_j > v(y,\,m),\,Y_l > v(y,\,m) \right]\\[2mm]
&&\leq (n-1)\,\sum_{i=1}^n  \sum_{j\neq l=1}^m\PP\left[ X_i >v(x,\,n),\;Y_j > v(y,\,m),\,Y_l > v(y,\,m)\right] \\[2mm]  
&&=o\left[\sum_{i=1}^n  \sum_{j\neq l=1}^m \left(\PP\left[ X_i >v(x,\,n),\;Y_j >v(y,\,m) \right]+\PP\left[ X_i >v(x,\,n),\;Y_l > v(y,\,m) \right]\right)\right],
\eeao
as $x \wedge y\to \infty$, thus from inclusion $\mathcal{C}^{(2)} \subsetneq \mathcal{D}^{(2)}$ we obtain
\beam \label{eq.KLP.3.10}
I_{22}(x,\,y) = o\left( \sum_{i=1}^n \sum_{j=1}^m \PP\left[ X_i >x,\;
Y_j >y \right] \right)\,,
\eeam
as $x \wedge y\to \infty$. Therefore from relations \eqref{eq.KLP.3.9} and \eqref{eq.KLP.3.10} in combination with \eqref{eq.KLP.3.8} we conclude
\beam \label{eq.KLP.3.11}
I_{2}(x,\,y) = o\left[S(x,\,y)\right]\,,
\eeam
as $x \wedge y\to \infty$, and then by \eqref{eq.KLP.3.7} and \eqref{eq.KLP.3.11} in combination with \eqref{eq.KLP.3.3} we get \eqref{eq.KLP.3.2}.

Now we estimate the lower bound in the form
\beam \label{eq.KLP.3.12}
\PP [ S_n >x,\;T_m>y] \geq [1-o(1)] \sum_{i=1}^n \sum_{j=1}^m \PP[ X_i >x,\;Y_j>y]\,,
\eeam  
as $x \wedge y\to \infty$. Let remind that $\mathcal{C}^{(2)} \subsetneq (\mathcal{D} \cap \mathcal{L})^{(2)}$, that means for any sequences $X_1,\,\ldots,\,X_n$ and $Y_1,\,\ldots,\,Y_m$, there exists some joint insensitivity function ${\bf a}$:
\beam \label{eq.KP.3.1*}
{\bf a}=(a_F,\,a_G):=\left(\bigwedge_{i=1}^n a_{F_i}\,,\; \bigwedge_{j=1}^m a_{G_j}\right)\,,
\eeam 
and this function represents an insensitivity function, for any distribution  pair $(F_i,\,G_j)$, for any $i=1,\,\ldots, \,n$ and $j=1,\,\ldots,\,m$. 

For such a function ${\bf a}$ we obtain
\beam \label{eq.KLP.3.13} 
\PP [ S_n >x,\,T_m>y] \geq  \PP\left[ S_n >x,\,T_m>y,\,\bigvee_{i=1}^n X_i >x+ a_F(x),\,\bigvee_{j=1}^mY_j>y+a_G(y) \right],
\eeam
Applying twice Bonferroni's inequality in \eqref{eq.KLP.3.13}, we find
\beam \label{eq.KLP.3.14} \notag
&&\PP [ S_n >x,\;T_m>y] \geq \sum_{i=1}^n \sum_{j=1}^m \PP\left[ S_n >x,\;T_m>y,\; X_i >x+ a_F(x),\;Y_j>y+a_G(y) \right]\\[2mm]
&&\qquad \qquad  -\sum_{i<k=1}^n \sum_{j=1}^m \PP\left[  X_i >x+ a_F(x),\;X_k>x+a_F(x),\; Y_j>y+a_G(y) \right]\\[2mm] \notag
&& \qquad \qquad  -\sum_{i=1}^n \sum_{j<l=1}^m  \PP\left[  X_i >x+ a_F(x),\; Y_j>y+a_G(y),\;Y_l>y+a_G(y) \right]\\[2mm] \notag
&&=:\sum_{k=1}^3 P_k(x,\,y) \,.
\eeam
By the assumption that  $X_1,\,\ldots,\,X_n$ and $Y_1,\,\ldots,\,Y_m$ are $GQAI$ and the insensitivity functions properties we conclude that 
\beam \label{eq.KLP.3.15} 
P_k(x,\,y)=o\left( \sum_{i=1}^n \sum_{j=1}^m \PP[ X_i >x,\,Y_j>y]\right)\,,   
\eeam
as $x \wedge y\to \infty$ for $k=2,\,3$. For $P_1(x,\,y)$ we write
\beam \label{eq.KLP.3.16} \notag
&&P_1(x,\,y) \geq  \sum_{i=1}^n \sum_{j=1}^m \PP\left[ X_i >x+ a_F(x),\;Y_j>y+a_G(y) \right]\\[2mm]
&& -\sum_{i\neq k=1}^n \sum_{j=1}^m \PP\left[  X_i >x+ a_F(x),\;X_k< -\dfrac{a_F(x)}n,\; Y_j>y+a_G(y) \right]\\[2mm] \notag
&& -\sum_{i=1}^n \sum_{j\neq l=1}^m \PP\left[  X_i >x+ a_F(x),\; Y_j>y+a_G(y),\;Y_l< -\dfrac{a_G(y)}m \right]\\[2mm] \notag
&&=:P_{11}(x,\,y) -P_{12}(x,\,y) -P_{13}(x,\,y)\,.
\eeam
Now we estimate $P_{12}(x,\,y)$
\beao
&&P_{12}(x,\,y) \lesssim \sum_{i\neq k=1}^n \sum_{j=1}^m \PP\left[ X_i > \dfrac{a_F(x)}n,\; Y_j>y,\;|X_k|>\dfrac{a_F(x)}n \right]\\[2mm] \notag
&& =o\left[ \sum_{i\neq k=1}^n \sum_{j=1}^m \left(\PP\left[  X_i >\dfrac{a_F(x)}n,\; Y_j>y \right] +\PP\left[ X_k > \dfrac{a_F(x)}n,\; Y_j>y\right] \right) \right]\\[2mm] \notag
&& =o\left( \sum_{i=1}^n  \sum_{j=1}^m \PP\left[  X_i >x,\; Y_j>y \right] \right) \,,
\eeao
as $x \wedge y\to \infty$, where we use the GQAI property and in last step we take into account the inclusion $\mathcal{C}^{(2)} \subsetneq \mathcal{D}^{(2)}$. By symmetry we find similarly 
\beam \label{eq.KLP.3.18} 
P_{13}(x,\,y) =o\left( \sum_{i=1}^n  \sum_{j=1}^m \PP\left[  X_i >x,\; Y_j>y \right] \right) \,,
\eeam
as $x \wedge y\to \infty$. Hence, from relations \eqref{eq.KLP.3.16} - \eqref{eq.KLP.3.18} and the fact that $\mathcal{C}^{(2)} \subsetneq (\mathcal{D}\cap \mathcal{L})^{(2)}$ in combination with the properties of function ${\bf a}$ we obtain
\beam \label{eq.KLP.3.19} 
P_{1}(x,\,y) \geq [1-o(1)] \sum_{i=1}^n  \sum_{j=1}^m \PP\left[  X_i >x,\; Y_j>y \right] \,,
\eeam
as $x \wedge y\to \infty$. So, by relation \eqref{eq.KLP.3.14} in combination with relations \eqref{eq.KLP.3.15} and \eqref{eq.KLP.3.19} we conclude that \eqref{eq.KLP.3.12} is true. Therefore
\beam \label{eq.KLP.3.20}
\PP[ S_n >x,\;T_m>y]\sim  S(x,\,y)\,,
\eeam 
as $x \wedge y\to \infty$. Now it remains to use \eqref{eq.KLP.3.20} and the inequality
\beao
\PP[ S_n >x,\;T_m>y]\leq \PP\left[ \bigvee_{i=1}^n S_i >x,\;\bigvee_{j=1}^m T_j>y\right] \leq \PP\left[  S_n^+ >x,\; T_m^+>y\right]\,,
\eeao
where 
\beao
S_n^+:=\sum_{i=1}^n X_i^+\,, \qquad T_m^+:=\sum_{j=1}^m Y_j^+\,,
\eeao 
to establish \eqref{eq.KLP.3.1}.
~\halmos

\bre \label{rem.KLP.3.1}
Comparing Theorem \ref{th.KLP.3.1} with relation \eqref{eq.KLP.1.8} we find the following differences. In relation \eqref{eq.KLP.1.8} the dependence structure is restricted as $GTAI \subsetneq GQAI$ (and in our Theorem we do not need the assumption of TAI in each sequence), while the distribution class is wider since $\mathcal{C}^{(2)} \subsetneq (\mathcal{D}\cap \mathcal{L})^{(2)}$. Furthermore relation \eqref{eq.KLP.1.8} contains also the asymptotic behavior of the maximums, which is missing from Theorem \ref{th.KLP.3.1}. Another difference is that in Theorem \ref{th.KLP.3.1} is not required the same multitude of summands, although in relation \eqref{eq.KLP.1.8} is not permitted, which can be helpful in risk theory, since it allows two different counting processes.
\ere

The two-dimensional risk models become more and more popular, because of the need of insurance companies to operate several parallel business lines, see for example in \cite{jiang:wang:chen:xu:2015}, \cite{li:yang:2015}, \cite{chen:cheng:zheng:2025}.

Next, we have a direct corollary for a two-dimensional discrete time risk model. The surplus processes are of the form 
\beao
U_1(k,\,x):=x- \sum_{i=1}^k X_i\,, \qquad  U_2(k,\,y):=y- \sum_{j=1}^k Y_j\,,
\eeao 
for any $k=1,\,\ldots,\,n$, where $x$ and $y$ are the initial surpluses in each of the two business lines, while $X_i$ and $Y_i$ are the net losses of each business line during the $i$-th period, namely total claims minus total premiums. In contrast to one-dimensional case, there exist several ways to define the ruin probability, see for example \cite{cheng:yu:2019}. Let choose the following definition of the ruin probability
\beao
\tau_1(x):= \inf \left\{k=1,\,\ldots,\,n\;:\; U_1(k,\,x)<0 \;|\; U_1(0)=x\right\}\,,\\[2mm]
\tau_2(y):= \inf \left\{k=1,\,\ldots,\,n\;:\; U_2(k,\,y)<0 \;|\; U_2(0)=y\right\}\,,
\eeao
and $\tau_{and}(x,\,y):=\max \{\tau_1(x),\,\tau_2(y)\}$, that indicates the first moment, when both surplus processes fall below zero, but not necessarily simultaneously. Hence we obtain
\beao
\psi_{and}(x,\,y;\,n):=\PP[ \tau_{and}(x,\,y) \leq n] = \PP\left[ \bigvee_{i=1}^n S_i >x,\;\bigvee_{j=1}^n T_j>y\right]\,.
\eeao 

Next we obtain an asymptotic expression of the ruin probability in finite time for $n \in \bbn$.

\bco \label{cor.KLP.3.1}
Under the conditions of Theorem \ref{th.KLP.3.1}, with $n \in \bbn$, then
\beao
\psi_{and}(x,\,y;\,n) \sim  \sum_{i=1}^n \sum_{j=1}^n \PP[ X_i >x,\;Y_j>y]\,,
\eeao 
as $x \wedge y\to \infty$.
\eco

\section {Two-dimensional Closure Properties}  \label{sec.KP.4}

Now we proceed to the closure properties of two-dimensional distributions, where some of them are very usefull in Section 5. Next lemma studies the closure properties of sum with respect to classes of regularly and consistently varying distributions under the dependence structure $QAI$. For the regular variation class the corresponding closure property under independent random variables can be found in \cite[p. 278]{feller:1971}, while for arbitrarily dependent random variables, under more strict conditions can be found in \cite[Lem. 4.4.2]{samorodnitsky:taqqu:1994} and \cite[Lem. 3.1]{leipus:surgailis:2007}. For the class of consistently varying distributions the closure property in the independence case was given by \cite[Th.2.2]{cai:tang:2004} for non-negative variables and by \cite[Lem. 3]{kizinevic:sprindys:siaulys:2016} for real valued random variables. Further, under more strict conditions in \cite[lem. 3.3(i)]{yang:yuen:liu:2018} we find closure property of sum under arbitrary dependence in class $\mathcal{C}$.

\ble \label{lem.KLP.4.1}
\begin{enumerate}
\item[(i)]
Let $X_1,\,X_2$ real valued, random variables,  with distributions $F_1,\,F_2 \in \mathcal{C}$ respectively, under $QAI$ dependence. Then
\beam \label{eq.KLP.4.1}
\bF_{X_1+X_2}(x) \sim  \bF_{1}(x) +\bF_{2}(x)\,,
\eeam 
as $\xto$, and $F_{X_1+X_2} \in \mathcal{C}$.
\item[(ii)]
Under the conditions of part $(i)$, with the only difference that $F_1,\,F_2 \in \mathcal{R}_{-\alpha}$, for some $\alpha >0$, then $F_{X_1+X_2} \in  \mathcal{R}_{-\alpha}$.
\end{enumerate} 
\ele

\pr~
\begin{enumerate}
\item[(i)]
In the case of $QAI$ relation \eqref{eq.KLP.4.1} follows from \cite[Th. 3.1]{chen:yuen:2009}, for $n=2$. Hence, for the closure property of $\mathcal{C}$ with respect to sum, we obtain for all the distributions
\beam \label{eq.KLP.4.2}
\lim_{z \uparrow 1} \limsup_{\xto} \dfrac {\bF_{X_1+X_2}(z\,x)}{\bF_{X_1+X_2}(x)} \geq 1\,,
\eeam
from the elementary inequalities 
\beam \label{eq.KLP.4.3}
\min \left\{ \dfrac ac\,,\;\dfrac bd \right\} \leq \dfrac {a+b}{c+d} \leq \max \left\{ \dfrac ac\,,\;\dfrac bd \right\} \,,
\eeam
for any constants $a,\,b,\,c,\,d >0$ and from \eqref{eq.KLP.4.1} we find
\beam \label{eq.KLP.4.4}
&&\lim_{z \uparrow 1} \limsup_{\xto} \dfrac {\bF_{X_1+X_2}(z\,x)}{\bF_{X_1+X_2}(x)} =\lim_{z \uparrow 1} \limsup_{\xto} \dfrac {\bF_{1}(z\,x)+\bF_{2}(z\,x)}{\bF_{1}(x)+\bF_{2}(x)} \\[2mm] \notag
&&\leq \lim_{z \uparrow 1} \limsup_{\xto}  \max \left\{ \dfrac {\bF_{1}(z\,x)}{\bF_{1}(x)}\,,\;\dfrac {\bF_{2}(z\,x)}{\bF_{2}(x)} \right\} =  \max_{i\in \{1,\,2\}} \left\{ \lim_{z \uparrow 1} \limsup_{\xto}\dfrac {\bF_{i}(z\,x)}{\bF_{i}(x)}\right\} =1\,,
\eeam
where in the last step we use the assumption $F_1,\,F_2 \in \mathcal{C}$. Hence by \eqref{eq.KLP.4.2} and \eqref{eq.KLP.4.4} we find $F_{X_1+X_2} \in \mathcal{C}$.
\item[(ii)]
From the relation $\mathcal{R}_{-\alpha}$ and by part $(i)$ we obtain \eqref{eq.KLP.4.1}, hence together with \eqref{eq.KLP.4.3} we conclude
\beam \label{eq.KLP.4.5}
\lim_{\xto} \dfrac {\bF_{X_1+X_2}(t\,x)}{\bF_{X_1+X_2}(x)} = \lim_{\xto}\dfrac {\bF_{1}(t\,x)+\bF_{2}(t\,x)}{\bF_{1}(x)+\bF_{2}(x)} \leq  \max_{i\in \{1,\,2\}} \left\{ \lim_{\xto}\dfrac {\bF_{i}(t\,x)}{\bF_{i}(x)}\right\} =t^{-\alpha}\,,
\eeam
for any  $t>0$, since $F_i \in \mathcal{R}_{-\alpha}$. Furthermore
\beam \label{eq.KLP.4.6}
\lim_{\xto} \dfrac {\bF_{X_1+X_2}(t\,x)}{\bF_{X_1+X_2}(x)} \geq  \min_{i\in \{1,\,2\}} \left\{ \lim_{\xto}\dfrac {\bF_{i}(t\,x)}{\bF_{i}(x)}\right\} =t^{-\alpha}\,,
\eeam
for any  $t>0$, so by \eqref{eq.KLP.4.5} and \eqref{eq.KLP.4.6} we find $F_{X_1+X_2} \in  \mathcal{R}_{-\alpha}$.~\halmos
\end{enumerate}

In the next result we find closure property of the distribution classes $\mathcal{C}^{(2)}$ and $\mathcal{R}_{(-\alpha_1,\,-\alpha_2)}^{(2)}$, with respect to sum under $GQAI$, with additional restriction that the particular summands are $QAI$. This way we have generalization of \cite[Cor. 4.1]{konstantinides:passalidis:2025} in case of $GQAI$ instead of $GTAI$ (and $QAI$ instead of TAI for the summands) and with real random variables.

\bth \label{th.KLP.4.1}
\begin{enumerate}
\item[(i)]
Let $X_1,\,X_2,\,Y_1,\,Y_2$ be real valued random variables,  with the  distributions $F_1,\,F_2,\,G_1,\,G_2$ respectively, under $GQAI$ structure. If $X_1,\,X_2$ are $QAI$ and $Y_1,\,Y_2$ are also $QAI$ with $(F_i,\,G_j) \in \mathcal{C}^{(2)}$ for any $i, \,j \in \{1,\,2\}$, then we find $\left(F_{X_1+X_2},\,G_{Y_1+Y_2}\right) \in \mathcal{C}^{(2)}$.
\item[(ii)]
Under the conditions of part $(i)$, with the only difference that $F_1,\,F_2 \in \mathcal{R}_{-\alpha_1}$ and $G_1,\,G_2 \in \mathcal{R}_{-\alpha_2}$, for some $\alpha_1,\,\alpha_2 >0$ where $(F_i,\,G_j) \in \mathcal{R}_{(-\alpha_1,\,-\alpha_2)}^{(2)}$, for any  $i, \,j \in \{1,\,2\}$, then $(F_{X_1+X_2},\;G_{Y_1+Y_2}) \in  \mathcal{R}_{(-\alpha_1,\,-\alpha_2)}^{(2)}$.
\end{enumerate}
\ethe

\pr~
\begin{enumerate}
\item[(i)]
From Theorem \ref{th.KLP.3.1}, for $n=m=2$ we get
\beam \label{eq.KLP.4.7}
\PP[X_1+X_2>x\,,\;Y_1+Y_2>y]\sim\sum_{i=1}^2 \sum_{j=1}^2 \PP[X_i>x\,,\;Y_j>y]\,,
\eeam
as $x \wedge y\to \infty$. For any two-dimensional distributions is true the lower bound
\beam \label{eq.KLP.4.8}
\lim_{{\bf z} \uparrow {\bf 1}} \limsup_{x \wedge y \to \infty} \dfrac {\PP[X_1+X_2>z_1\,x\,,\;Y_1+Y_2>z_2\,y]}{\PP[X_1+X_2>x\,,\;Y_1+Y_2>y]} \geq 1 \,,
\eeam
while by \eqref{eq.KLP.4.7}, applying the upper inequality in \eqref{eq.KLP.4.3} thrice, we obtain
\beam \label{eq.KLP.4.9} \notag
&&\lim_{{\bf z} \uparrow {\bf 1}} \limsup_{x \wedge y \to \infty} \dfrac {\PP[X_1+X_2>z_1\,x\,,\;Y_1+Y_2>z_2\,y]}{\PP[X_1+X_2>x\,,\;Y_1+Y_2>y]} \\[2mm]
&&=
\lim_{{\bf z} \uparrow {\bf 1}} \limsup_{x \wedge y \to \infty} \dfrac {\sum_{i=1}^2 \sum_{j=1}^2 \PP[X_i>z_1\,x\,,\;Y_j>z_2\,y]}{\sum_{i=1}^2 \sum_{j=1}^2 \PP[X_i>x\,,\;Y_j>y]} \\[2mm] \notag
&&\leq \max_{i,\,j\in \{1,\,2\}} \left\{ \lim_{{\bf z} \uparrow {\bf 1}} \limsup_{x \wedge y \to \infty} \dfrac {\PP[X_i>z_1\,x\,,\;Y_j>z_2\,y]}{\PP[X_i>x\,,\;Y_j>y]}\right\}=1\,,
\eeam
where in the last step we used that $(F_i,\,G_j) \in \mathcal{C}^{(2)}$, for any  $i, \,j \in \{1,\,2\}$. Hence, by \eqref{eq.KLP.4.8} and \eqref{eq.KLP.4.9} we find \eqref{eq.KLP.2.5}. Furthermore, by Lemma \ref{lem.KLP.4.1}(i) we obtain $F_{X_1+X_2}  \in \mathcal{C}$ and $G_{Y_1+Y_2}  \in \mathcal{C}$, which together with \eqref{eq.KLP.2.5} gives $(F_{X_1+X_2},\,G_{Y_1+Y_2}) \in \mathcal{C}^{(2)}$.
\item[(ii)]
From $\mathcal{R}_{(-\alpha_1,\,-\alpha_2)}^{(2)} \subsetneq  \mathcal{C}^{(2)}$
, using Theorem \ref{th.KLP.3.1}, with $n=m=2$, we reach to \eqref{eq.KLP.4.7}, and from this, through the application of upper inequality in \eqref{eq.KLP.4.3} thrice, we finally get
\beam \label{eq.KLP.4.10} \notag
&& \lim_{x \wedge y \to \infty} \dfrac {\PP[X_1+X_2>t_1\,x\,,\;Y_1+Y_2>t_2\,y]}{\PP[X_1+X_2>x\,,\;Y_1+Y_2>y]} =
 \lim_{x \wedge y \to \infty} \dfrac {\sum_{i=1}^2 \sum_{j=1}^2 \PP[X_i>t_1\,x\,,\;Y_j>t_2\,y]}{\sum_{i=1}^2 \sum_{j=1}^2 \PP[X_i>x\,,\;Y_j>y]} \\[2mm] 
&&\qquad \qquad \qquad \leq \max_{i,\,j\in \{1,\,2\}} \left\{\limsup_{x \wedge y \to \infty} \dfrac {\PP[X_i>t_1\,x\,,\;Y_j>t_2\,y]}{\PP[X_i>x\,,\;Y_j>y]}\right\}=t_1^{-\alpha_1}\,t_2^{-\alpha_2}\,,
\eeam
for any $t_1,\,t_2 > 0$, where in last step was used that $(F_i,\,G_j) \in \mathcal{R}_{(-\alpha_1,\,-\alpha_2)}^{(2)}$, for any  $i, \,j \in \{1,\,2\}$, and with similar way, by application of lower inequality in \eqref{eq.KLP.4.3}, we have
\beam \label{eq.KLP.4.11} \notag
\lim_{x \wedge y \to \infty} \dfrac {\PP[X_1+X_2>t_1\,x\,,\;Y_1+Y_2>t_2\,y]}{\PP[X_1+X_2>x\,,\;Y_1+Y_2>y]} &\geq&  \min_{i,\,j\in \{1,\,2\}} \left\{\lim_{x \wedge y \to \infty} \dfrac {\PP[X_i>t_1\,x\,,\;Y_j>t_2\,y]}{\PP[X_i>x\,,\;Y_j>y]}\right\} \\[2mm]
&=&t_1^{-\alpha_1}\,t_2^{-\alpha_2}\,,
\eeam
thus by \eqref{eq.KLP.4.10} and  \eqref{eq.KLP.4.11} we have  \eqref{eq.KLP.2.6}. Next, by Lemma \ref{lem.KLP.4.1}(ii) we obtain the inclusions $F_{X_1+X_2} \in \mathcal{R}_{-\alpha_1}$ and $F_{Y_1+Y_2} \in \mathcal{R}_{-\alpha_2}$, that together with \eqref{eq.KLP.2.6} finally gives $(F_{X_1+X_2},\;G_{Y_1+Y_2}) \in  \mathcal{R}_{(-\alpha_1,\,-\alpha_2)}^{(2)}$.~\halmos
\end{enumerate}

Now we study the closure properties of distribution classes with respect to product convolution in two dimensions. This work can help the extension of Theorem \ref{th.KLP.3.1} to the direction of random weighted sums. 

To show the closure property with respect to convolution product in the distributions classes $\mathcal{D}^{(2)}$, $\mathcal{C}^{(2)}$, $\mathcal{L}^{(2)}$ we employ the following assumption.

\begin{assumption} \label{ass.KP.A}
Let $b\,:\,[0,\,\infty) \longrightarrow (0,\,\infty)$ be a function, such that
$b(x) \to \infty$, $b(x)=o(x)$, as $\xto$ and
\beam \label{eq.KLP.4.14} 
\PP[\Theta>b(x)]=o(\PP[\Theta\,X>x\,,\;\Delta\,Y>y])\,,
\eeam
as $x \wedge y \to \infty$, and  $c\,:\,[0,\,\infty) \longrightarrow (0,\,\infty)$ be a function, such that
$c(y) \to \infty$, $c(y)=o(y)$, as $\yto$ and
\beam \label{eq.KLP.4.15} 
\PP[\Delta>c(y)]=o(\PP[\Theta\,X>x\,,\;\Delta\,Y>y])\,,
\eeam
as $x \wedge y \to \infty$.
\end{assumption}

\bre \label{rem.KLP.4.2}
It is easy to see that if the $\Theta$ and $\Delta$ have distributions with  upper bounded support, because $X,\,Y$ are heavy tailed in this paper, hence have unbounded from above supports, then relations \eqref{eq.KLP.4.14} and \eqref{eq.KLP.4.15} are true directly. We observe that by relation \eqref{eq.KLP.4.14} follows 
\beam \label{eq.KLP.4.16} 
\PP[\Theta>b(x)]=o(\PP[\Theta\,X>x])\,,
\eeam
as $\xto$, and by relation \eqref{eq.KLP.4.15} follows 
\beam \label{eq.KLP.4.17} 
\PP[\Delta>c(y)]=o(\PP[\Delta\,Y>y])\,,
\eeam
as $\yto$. Therefore, from \eqref{eq.KLP.4.16} and \eqref{eq.KLP.4.17}, in combination with the definition of functions $b(x)$ and $c(y)$, by \cite[Lem. 3.2]{tang:2006} follows that if $\Theta$ and $\Delta$ have distributions with unbounded supports, then
\beam \label{eq.KLP.4.18} 
\PP[\Theta>u\,x]=o(\PP[\Theta\,X>x])\,, \qquad \PP[\Delta>u\,y]=o(\PP[\Delta\,Y>y])\,,
\eeam
as $\xto$ and as $\yto$ respectively, for any $u>0$. In case $\Theta$ and $\Delta$ have distributions with upper bounded supports, relation \eqref{eq.KLP.4.18} is true, when $X$ and  $Y$ are heavy-tailed, or in general if they have unbounded from above supports.
\ere

\bth \label{th.KLP.4.2}
Let $(X,\,Y)$ be a random pair with the  distributions $F,\,G$ respectively, $(\Theta,\,\Delta)$ be a random pair, independent of $(X,\,Y)$, with non-negative, non-degenerate to zero marginal distributions, and Assumption \ref{ass.KP.A} hold. The following are true
\begin{enumerate}
\item[(i)]
If  $(X,\,Y) \in \mathcal{D}^{(2)}$, then  $(\Theta\,X,\,\Delta\,Y) \in \mathcal{D}^{(2)}$.
\item[(ii)]
If  $(X,\,Y) \in \mathcal{L}^{(2)}$, then  $(\Theta\,X,\,\Delta\,Y) \in \mathcal{L}^{(2)}$.
\item[(iii)]
If  $(X,\,Y) \in \mathcal{C}^{(2)}$, then  $(\Theta\,X,\,\Delta\,Y) \in \mathcal{C}^{(2)}$.
\item[(iv)]
If  $(X,\,Y) \in (\mathcal{D}\cap \mathcal{L})^{(2)}$, then  $(\Theta\,X,\,\Delta\,Y) \in (\mathcal{D}\cap \mathcal{L})^{(2)}$.
\end{enumerate}
\ethe

\pr~
\begin{enumerate}
\item[(i)]
Let  ${\bf b}=(b_1,\,b_2) \in (0,\,1)^{2}$, then  
\beam \label{eq.KLP.4.19}  \notag
&&\PP[\Theta\,X>b_1\,x,\;\Delta\,Y>b_2\,y]=\int_0^{\infty} \int_0^{\infty} \PP\left[X>\dfrac{b_1\,x}t,\;Y>\dfrac{b_2\,y}s\right]\,\PP[\Theta \in dt,\;\Delta \in ds]\\[2mm] \notag
&&=\left(\int_0^{b(x)}+\int_{b(x)}^{\infty}\right) \left(\int_0^{c(y)}+\int_{c(y)}^{\infty}\right) \PP\left[X>\dfrac{b_1\,x}t,\;Y>\dfrac{b_2\,y}s\right]\,\PP[\Theta \in dt,\;\Delta \in ds]\\[2mm]
&&=:I_{11}(x,\,y)+I_{12}(x,\,y)+I_{21}(x,\,y)+I_{22}(x,\,y)\,. 
\eeam
But from Assumption  \ref{ass.KP.A} we obtain
\beam \label{eq.KLP.4.20} 
I_{22}(x,\,y)&=&\int_{b(x)}^{\infty} \int_{c(y)}^{\infty} \PP\left[X>\dfrac{b_1\,x}t,\;Y>\dfrac{b_2\,y}s\right]\,\PP[\Theta \in dt,\;\Delta \in ds]\\[2mm] \notag
&\leq& \PP[\Theta>b(x),\;\Delta>c(y)]\leq \PP[\Theta>b(x)]=o(\PP[\Theta\,X>x,\;\Delta\,Y>y])\,, 
\eeam
as $x \wedge y \to \infty$. For $I_{12}(x,\,y)$  we find
\beam  \label{eq.KLP.4.21} 
I_{12}(x,\,y)\leq \PP[\Theta\leq b(x),\;\Delta>c(y)]\leq \PP[\Delta>c(y)]=o(\PP[\Theta\,X>x,\;\Delta\,Y>y])\,, 
\eeam
as $x \wedge y \to \infty$. Through symmetry, we obtain similarly
\beam \label{eq.KLP.4.22} 
I_{21}(x,\,y)\leq \PP[\Theta>b(x),\;\Delta\leq c(y)]\leq \PP[\Theta>b(x)]=o(\PP[\Theta\,X>x,\;\Delta\,Y>y])\,, 
\eeam
as $x \wedge y \to \infty$. Hence, by \eqref{eq.KLP.4.19} together with \eqref{eq.KLP.4.20} - \eqref{eq.KLP.4.22} follows
\beam \label{eq.KLP.4.23} \notag
&&\limsup_{x \wedge y \to \infty}\dfrac{\PP[\Theta\,X>b_1\,x,\;\Delta\,Y>b_2\,y]}{\PP[\Theta\,X>x,\;\Delta\,Y>y]}\\[2mm] \notag
&&\leq \limsup_{x \wedge y \to \infty}\dfrac{\int_0^{b(x)} \int_0^{c(y)} \PP\left[X>\dfrac{b_1\,x}t,\;Y>\dfrac{b_2\,y}s\right]\,\PP[\Theta \in dt,\;\Delta \in ds]}{\int_0^{b(x)} \int_0^{c(y)} \PP\left[X>\dfrac{x}t,\;Y>\dfrac{y}s\right]\,\PP[\Theta \in dt,\;\Delta \in ds]} + o(1)\\[2mm]
&&\leq \limsup_{x \wedge y \to \infty} \sup_{t \in (0,\,b(x)],\,s \in (0,\,c(y)]}\dfrac{\PP\left[X>\dfrac{b_1\,x}t,\;Y>\dfrac{b_2\,y}s\right]}{ \PP\left[X>\dfrac{x}t,\;Y>\dfrac{y}s\right]}+o(1) \\[2mm] \notag
&&= \limsup_{x \wedge y \to \infty}\dfrac{\PP\left[X>b_1\,x,\;Y>b_2\,y\right]}{ \PP\left[X>x,\;Y>y \right]}< \infty\,, 
\eeam
where in the last step used that $(X,\,Y) \in \mathcal{D}^{(2)}$. Next, since $F\in \mathcal{D}$ and $G \in \mathcal{D}$, through \cite[Th. 3.3(ii)]{cline:samorodnitsky:1994} or \cite[Prop. 5.4(i)]{leipus:siaulys:konstantinides:2023} we conclude that $\Theta\,X \in \mathcal{D}$ and $\Delta\,Y \in \mathcal{D}$, which in combination with \eqref{eq.KLP.4.23} provided $(\Theta\,X,\;\Delta\,Y) \in \mathcal{D}^{(2)}$.
\item[(ii)]
In case $(X,\,Y) \in \mathcal{L}^{(2)}$, for $a_1,\,a_2 >0$ is well known that for any two-dimensional distribution
\beao 
\liminf_{x \wedge y\to \infty}\dfrac{\PP[\Theta\,X>x-a_1,\;\Delta\,Y>y-a_2]}{\PP[\Theta\,X>x,\;\Delta\,Y>y]}\geq 1\,. 
\eeao
On the other hand side we find
\beam \label{eq.KLP.4.28} 
&& \PP[\Theta\,X>x-a_1,\;\Delta\,Y>y-a_2]\\[2mm] \notag
&& =\left(\int_0^{b(x)}+\int_{b(x)}^{\infty}\right) \left(\int_0^{c(y)}+\int_{c(y)}^{\infty}\right) \PP\left[X>\dfrac{x-a_1}t,\;Y>\dfrac{y-a_2}s\right]\,\PP[\Theta \in dt,\;\Delta \in ds]\\[2mm] \notag
&&=: K_{11}(x,\,y)+K_{12}(x,\,y)+K_{21}(x,\,y)+K_{22}(x,\,y)\,. 
\eeam
Here also, similarly to  \eqref{eq.KLP.4.20} - \eqref{eq.KLP.4.22}, we find that the $K_{12}(x,\,y)$, $K_{21}(x,\,y)$, $K_{22}(x,\,y)$ are of $o(\PP[\Theta\,X>x,\;\Delta\,Y>y])$ order of magnitude, as $x \wedge y\to \infty$. Thus, together with  \eqref{eq.KLP.4.28} we conclude
\beam \label{eq.KLP.4.29} \notag
&&\limsup_{x \wedge y \to \infty}\dfrac{\PP[\Theta\,X>x-a_1,\;\Delta\,Y>y-a_2]}{\PP[\Theta\,X>x,\;\Delta\,Y>y]}\\[2mm]  \notag
&& \leq \limsup_{x \wedge y \to \infty}\dfrac{\int_0^{b(x)}\int_0^{c(y)} \PP\left[X>\dfrac{x-a_1}t,\;Y>\dfrac{y-a_2}s\right]\,\PP[\Theta \in dt,\;\Delta \in ds]}{\int_0^{b(x)}\int_0^{c(y)} \PP\left[X>\dfrac{x}t,\;Y>\dfrac{y}s\right]\,\PP[\Theta \in dt,\;\Delta \in ds]}+o(1)\\[2mm]
&& \leq \limsup_{x \wedge y \to \infty}\sup_{t \in (0,\,b(x)],\,s \in (0,\,c(y)]}\dfrac{\PP\left[X>\dfrac{x-a_1}t,\;Y>\dfrac{y-a_2}s\right]}{ \PP\left[X>\dfrac{x}t,\;Y>\dfrac{y}s\right]}\\[2mm] \notag
&&=\limsup_{x \wedge y \to \infty} \dfrac{\PP[X>x-a_1,\;Y>y-a_2]}{\PP[X>x,\;Y>y]}= 1\,, 
\eeam
where in the last step we take into account that $(X,\,Y) \in \mathcal{L}^{(2)}$. By Assumption \ref{ass.KP.A} and Remark \ref{rem.KLP.4.2}, through \cite[Cor. 5.1]{leipus:siaulys:konstantinides:2023}  we obtain $\Theta\,X \in \mathcal{L}$, $\Delta\,Y \in \mathcal{L}$, which together with \eqref{eq.KLP.4.29} provides $(\Theta\,X,\,\Delta\,Y) \in \mathcal{L}^{(2)}$.
\item[(iii)]
For any two-dimensional distribution holds the inequality
\beam \label{eq.KLP.4.24} 
\lim_{{\bf z} \uparrow {\bf 1}} \limsup_{x \wedge y\to \infty}\dfrac{\PP[\Theta\,X>z_1\,x,\;\Delta\,Y>z_2\,y]}{\PP[\Theta\,X>x,\;\Delta\,Y>y]}\geq 1\,, 
\eeam
from the other hand side we obtain
\beam \label{eq.KLP.4.25} 
&& \PP[\Theta\,X>z_1\,x,\;\Delta\,Y>z_2\,y] \\[2mm] \notag
&&=\left(\int_0^{b(x)}+\int_{b(x)}^{\infty}\right) \left(\int_0^{c(y)}+\int_{c(y)}^{\infty}\right) \PP\left[X>\dfrac{z_1\,x}t,\;Y>\dfrac{z_2\,y}s\right]\,\PP[\Theta \in dt,\;\Delta \in ds]\\[2mm] \notag
&&=: J_{11}(x,\,y)+J_{12}(x,\,y)+J_{21}(x,\,y)+J_{22}(x,\,y)\,. 
\eeam
Similarly to \eqref{eq.KLP.4.20} - \eqref{eq.KLP.4.22}, we find that the $J_{12}(x,\,y)$, $J_{21}(x,\,y)$, $J_{22}(x,\,y)$ are of $
o(\PP[\Theta\,X>x,\;\Delta\,Y>y])$ order of magnitude, as $x \wedge y \to \infty$. Therefore, together with relation \eqref{eq.KLP.4.25} we obtain
\beam \label{eq.KLP.4.26}  \notag
&& \dfrac{\PP[\Theta\,X>z_1\,x,\,\Delta\,Y>z_2\,y]}{\PP[\Theta\,X>x,\,\Delta\,Y>y]}\lesssim \dfrac{\int_0^{b(x)}\int_0^{c(y)} \PP\left[X>\dfrac{z_1\,x}t,\,Y>\dfrac{z_2\,y}s\right]\,\PP[\Theta \in dt,\,\Delta \in ds]}{\int_0^{b(x)}\int_0^{c(y)} \PP\left[X>\dfrac{x}t,\;Y>\dfrac{y}s\right]\,\PP[\Theta \in dt,\;\Delta \in ds]}\\[2mm] \notag
&&+o(1)\lesssim  \sup_{t \in (0,\,b(x)],\,s \in (0,\,c(y)]}\dfrac{\PP\left[X>\dfrac{z_1\,x}t,\;Y>\dfrac{z_2\,y}s\right]}{ \PP\left[X>\dfrac{x}t,\;Y>\dfrac{y}s\right]}=\dfrac{\PP[X>z_1\,x,\;Y>z_2\,y]}{\PP[X>x,\;Y>y]}\rightarrow 1\,, \\
\eeam
as ${\bf z} \uparrow {\bf 1}$ and $x \wedge y\to \infty$, where in last step was used the condition $(X,\,Y) \in \mathcal{C}^{(2)}$. Now taking into consideration Assumption  \ref{ass.KP.A} and $X \in \mathcal{C}$, $Y \in \mathcal{C}$, through Remark \ref{rem.KLP.4.2}, because of \eqref{eq.KLP.4.24}, via \cite[Th. 3.4(ii)]{cline:samorodnitsky:1994}, or \cite[Prop. 5.3(ii)]{leipus:siaulys:konstantinides:2023} we obtain $\Theta \,X \in \mathcal{C}$, $\Delta \,Y \in \mathcal{C}$, which together with \eqref{eq.KLP.4.24} and \eqref{eq.KLP.4.26} provides $(\Theta\,X,\,\Delta\,Y) \in \mathcal{C}^{(2)}$.
\item[(iv)]
This case follows from parts $(i)$ and $(ii)$.~\halmos
\end{enumerate}

In part $(ii)$ of Theorem \ref{th.KLP.4.2} we find a generalization of \cite[Lem. 6.2]{konstantinides:passalidis:2025}, since now we have real-valued random variables $X$ and $Y$. 
In next lemma we find an extension of Breiman's Theorem in case of distribution class $\mathcal{R}_{(-\alpha_1,\,-\alpha_2)}^{(2)}$. For the class of regular variation in one dimensional case, see \cite{breiman:1965}, \cite{konstantinides:mikosch:2005}, \cite{denisov:zwart:2007},  while for the class of multivariate regular variation, see \cite{basrak:davis:mikosch:2002b}, \cite{fougeres:mercadier:2012}.

\ble \label{lem.KLP.4.2}
Let $(X,\,Y)$ be a random pair with the  distributions $F \in \mathcal{R}_{-\alpha_{1}}$, $G \in \mathcal{R}_{-\alpha_{2}}$ respectively, such that $(F,\;G) \in  \mathcal{R}_{(-\alpha_1,\,-\alpha_2)}^{(2)}$, with $0<\alpha_1,\,\alpha_2 < \infty$. Let $(\Theta,\,\Delta)$ be a random pair, independent of $(X,\,Y)$, with non-negative and non-degenerate to zero distributions, such that $\E\left[\Theta^{\alpha_1+\vep}\,\Delta^{\alpha_2+\vep}\right] < \infty$, for some $\vep >0$. Under the Assumption \ref{ass.KP.A} we have 
\beam \label{eq.KLP.4.12} 
\PP[\Theta\,X>x\,,\;\Delta\,Y>y] \sim \E\left[\Theta^{\alpha_1}\,\Delta^{\alpha_2}\right]\,\PP[X>x\,,\;Y>y] \,,
\eeam
as $x \wedge y\to \infty$. Furthermore  $(\Theta\,X,\;\Delta\,Y) \in  \mathcal{R}_{(-\alpha_1,\,-\alpha_2)}^{(2)}$.
\ele

\pr~
We have
\beao
&&{\PP[\Theta\,X> x\,,\;\Delta\,Y>y]} = \int_0^{\infty}\int_0^{\infty}  {\PP\left[X> \dfrac xt\,,\;Y> \dfrac ys\right]}\,\PP[\Theta \in dt\,,\;\Delta \in ds] \\[2mm]
&&={\left(\int_0^{b(x)}+\int_{b(x)}^{\infty}\right) \left(\int_0^{c(y)}+\int_{c(y)}^{\infty}\right) \PP\left[X>\dfrac{x}t,\;Y>\dfrac{y}s\right]\,\PP[\Theta \in dt,\;\Delta \in ds]} \\[2mm] 
&&=: \Lambda_{11}(x,y)+\Lambda_{1,2}(x,y)+\Lambda_{2,1}(x,y)+\Lambda_{2,2}(x,y)\,.
\eeao
By Assumption \ref{ass.KP.A} and similarly to \eqref{eq.KLP.4.20} - \eqref{eq.KLP.4.22}, we find that the $\Lambda_{12}(x,\,y)$, $\Lambda_{21}(x,\,y)$, $\Lambda_{22}(x,\,y)$ are of 
\beao
o(\PP[\Theta\,X>x,\;\Delta\,Y>y])\,,
\eeao 
order of magnitude, as $x \wedge y\to \infty$. Hence, by definitions of functions $b,c$, and since for some $\vep >0$ holds $\E\left[\Theta^{\alpha_1+\vep}\,\Delta^{\alpha_2+\vep}\right] < \infty$, from dominated convergence theorem and by $(F,\;G) \in  \mathcal{R}_{(-\alpha_1,\,-\alpha_2)}^{(2)}$ we have 
\beao
&& {\PP[\Theta\,X> x\,,\;\Delta\,Y>y]} \sim [1- o(1)]\,\int_0^{b(x)}\int_0^{c(y)}  {\PP\left[X> \dfrac xt\,,\;Y> \dfrac ys\right]}\,\PP[\Theta \in dt\,,\;\Delta \in ds] \\[2mm] 
&&\sim\int_0^{b(x)}\int_0^{c(y)} t^{\alpha_1}\,s^{\alpha_2} {\PP\left[X> x\,,\;Y>  y \right]}\,\PP[\Theta \in dt\,,\;\Delta \in ds]\\[2mm]
&&= E\left[\Theta^{\alpha_1}\Delta^{\alpha_2}{\bf 1}_{\{\Theta\leq b(x),\Delta\leq c(y)\}}\right] \PP\left[X> x\,,\;Y>  y \right]\,,
\eeao
as $x \wedge y\to \infty$, by the definition of functions $b,c$, namely tends to infinity, the indicator function tends to unity, as a result we have the relation \eqref{eq.KLP.4.12}.
Further, for any $t_1,t_2>0$ by \eqref{eq.KLP.4.12} we have
\beao
\dfrac{\PP[\Theta\,X>t_1\,x,\,\Delta\,Y>t_2\,y]}{\PP[\Theta\,X>x,\,\Delta\,Y>y]}= \dfrac{\E\left[\Theta^{\alpha_1}\,\Delta^{\alpha_2}\right]\PP[X>t_1\,x,\,Y>t_2\,y]}{\E\left[\Theta^{\alpha_1}\,\Delta^{\alpha_2}\right]\PP[X>x,\,Y>y]}\rightarrow t_{1}^{-\alpha_1}\,t_{2}^{-\alpha_2}\,,
\eeao
as $x \wedge y\to \infty$, where in the last step we used the fact that $(X\;Y) \in  \mathcal{R}_{(-\alpha_1,\,-\alpha_2)}^{(2)}$ By the assumption $\E\left[\Theta^{\alpha_1+\vep}\,\Delta^{\alpha_2+\vep}\right] < \infty$, for some $\vep >0$, we have that $\E\left[\Theta^{\alpha_1+\vep}\right] < \infty$, and $\E\left[\Delta^{\alpha_2+\vep}\right] < \infty$, for some $\vep >0$, and by Breiman's Theorem (see for example Proposition 5.2 of \cite{leipus:siaulys:konstantinides:2023}), we have that $\Theta X\in\mathcal{R}_{-\alpha_{1}}$ and  $\Delta Y\in\mathcal{R}_{-\alpha_{2}}$, which in combination with the last equation provides $
(\Theta\,X,\;\Delta\,Y) \in  \mathcal{R}_{(-\alpha_1,\,-\alpha_2)}^{(2)}$.~\halmos

\bre \label{rem.KLP.4.1}
The previous result, except closure property, provides property  \eqref{eq.KLP.4.12} as well, which can be helpful to direct asymptotic expression of the joint asymptotic behavior of the randomly weighted sums, in the frame of class $\mathcal{R}_{(-\alpha_1,\,-\alpha_2)}^{(2)}$. We observe that the random weights $\Theta,\, \Delta$ are arbitrarily dependent, that can play crucial role in actuarial applications, since they represent discount factors. 
\ere

\section{Randomly Weighted Sums}  \label{sec.KP.5}

In this section we study the joint behavior of randomly weighted sums. There are several papers about formula \eqref{eq.KLP.1.7}, see for example \cite{li:2018b}, \cite{shen:ge:fu:2020}, \cite{shen:du:2023}, \cite{yang:chen:yuen:2024}. This joint behavior can be applied on several areas of actuarial science and financial mathematics. Indeed, in risk theory the main components $X_i$ and $Y_j$ represent gains or losses in $i$-th  and $j$-th period, while the random weights $\Theta_i$ and $\Delta_j$ represent discount factors, that are allowed to be degenerated to some positive number, see for example \cite{wang:su:yang:2024}. In credit risk applications, the random variables $X_i$ and $Y_j$ represent the rate of default of the $i$-th and $j$-th obligor and the random weights $\Theta_i$ and $\Delta_j$ represent Bernoulli random variables with values zero and unity, where the unity reflect the case of default.

Before the main result, we need a lemma, which examines the closure property for the dependencies $GTAI$ and $GQAI$ with respect to the convolution product under the distribution class $\mathcal{D}^{(2)}$ for the main variables, see some results for $TAI$ and $QAI$ in \cite[Th. 2.2]{li:2013}. We should mention that from here and later, when we say that the pair $(\Theta_i,\,\Delta_j)$ satisfies Assumption \ref{ass.KP.A}, we have in mind that $b_i(x_i) \to \infty$, $b_i(x_i)=o(x_i)$, as $x_i \to \infty$, for any $i=1,\,\ldots,\,n$, $c_j(y_j) \to \infty$, $c_j(y_j)=o(y_j)$, as $y_j \to \infty$, for any $j=1,\,\ldots,\,m$ and the relations \eqref{eq.KLP.4.14} and \eqref{eq.KLP.4.15} hold for any pair $(\Theta_i\,X_i\,,\;\Delta_j\,Y_j)$, for any  $i=1,\,\ldots,\,n$,  $j=1,\,\ldots,\,m$. 

\ble \label{lem.KLP.5.1}
Let $\{X_1,\,\ldots,\,X_n\}$, $\{Y_1,\,\ldots,\,Y_m\}$ be real valued, random variables with corresponding distributions $F_1,\,\ldots,\,F_n,\,G_1,\,\ldots,\,G_m$, which are $GTAI$ (or $GQAI$) and for any $i=1,\,\ldots,\,n$, $j=1,\,\ldots,\,m$ hold $(X_i,\,Y_j) \in \mathcal{D}^{(2)}$. Let $\Theta_1,\,\ldots,\,\Theta_n,\,\Delta_1,\,\ldots,\,\Delta_m$ be non-negative, non-degenerate to zero random variables with $(\Theta_i,\,\Delta_j)$ satisfying Assumption \ref{ass.KP.A}. If additionally we assume  $\Theta_1,\,\ldots,\,\Theta_n,\,\Delta_1,\,\ldots,\,\Delta_m$ are independent of $X_1,\,\ldots,\,X_n,\,Y_1,\,\ldots,\,Y_m$, then the products $\Theta_1\,X_1,\,\ldots,\,\Theta_n\,X_n,\,\Delta_1\,Y_1,\,\ldots,\,\Delta_m\,Y_m$ are $GTAI$, (or $GQAI$, respectively).
\ele

\pr~
Let begin with the $GTAI$ case. We denote the maximum $\widehat{\Theta}_{ik}:=\Theta_i \vee \Theta_k$, for $i\neq k \in \{1,\,\ldots,\,n\}$ and by Assumption \ref{ass.KP.A}, and $\widehat{b}(x_i \wedge x_k):=\max \{b_i(x_i \wedge x_k),\,b_k(x_i \wedge x_k)\}$, where the functions $b_i$ and $b_k$ are defined according to Assumption   \ref{ass.KP.A}  for the random weights $\Theta_i$ and $\Theta_k$, respectively, we find that $\widehat{b}(x_i \wedge x_k) \to \infty$ and $\widehat{b}(x_i \wedge x_k)=o(x_i \wedge x_k)$, as $x_i \wedge x_k \to \infty$ and further for some $j=1,\,\ldots,\,m$
\beam \label{eq.KLP.5.1} 
\PP\left[\widehat{\Theta}_{ik}>\widehat{b}(x_i \wedge x_k)\right]=o(\PP[\Theta_k\,X_k>x,\;\Delta_j\,Y_j>y])\,, 
\eeam
as $x_i \wedge x_k \to \infty$. Indeed, if $\widehat{\Theta}_{ik}=\Theta_k$, then by Assumption   \ref{ass.KP.A} relation \eqref{eq.KLP.5.1} is obvious. If $\widehat{\Theta}_{ik} =\Theta_i$, we multiply and divide the left member of \eqref{eq.KLP.5.1}  with $\PP[\Theta_i\,X_i >x_i,\,\Delta_j\,Y_j>y_j]$ to find that
\beao
\dfrac{\PP\left[\Theta_{i}>\widehat{b}(x_i \wedge x_k)\right]}{\PP[\Theta_i\,X_i >x_i,\,\Delta_j\,Y_j>y_j]}\,,
\eeao
tends to zero, while
\beao
\dfrac{\PP[\Theta_i\,X_i >x_i,\,\Delta_j\,Y_j>y_j]}{\PP[\Theta_k\,X_k >x_k,\,\Delta_j\,Y_j>y_j]} \leq \dfrac 1{\PP[\Theta_k\,X_k >x_k,\,\Delta_j\,Y_j>y_j]}<\infty\,,
\eeao 
where the last inequality holds because $X,Y$ have distributions with unbounded supports.

Hence, for a function $c_j$, defined in Assumption \ref{ass.KP.A} for the random variable $Y_j$, for any $j=1,\,\ldots,\,m$ and for $x_i,\,x_k,\,y_j >0$ we obtain
\beam \label{eq.KLP.5.2} \notag
&&\PP\left[|\Theta_i\,X_i|>x_i,\;\Theta_k\,X_k>x_k,\; \Delta_j\,Y_j>y_j\right]\leq \PP\left[\widehat{\Theta}_{ik}\,|X_i|>x_i,\; \widehat{\Theta}_{ik}\,X_k>x_k,\;\Delta_j\,Y_j>y_j\right]\\[2mm] 
&&=\int_0^{\infty} \int_0^{\infty}   \PP\left[|X_i|>\dfrac{x_i}t,\; X_k>\dfrac {x_k}t,\;Y_j>\dfrac {y_j}s\right]\,\PP\left[\widehat{\Theta}_{ik}\in dt,\;\Delta_j\in ds\right]\\[2mm] \notag
&&=\left(\int_0^{\widehat{b}(x_i \wedge x_k)}+\int_{\widehat{b}(x_i \wedge x_k)}^{\infty}\right) \left(\int_0^{c_j(y_j)}+\int_{c_j(y_j)}^{\infty}\right)  \PP\left[|X_i|>\dfrac{x_i}t,\; X_k>\dfrac {x_k}t,\;Y_j>\dfrac {y_j}s\right]\\[2mm] \notag
&& \qquad \qquad \times \PP\left[\widehat{\Theta}_{ik}\in dt,\;\Delta_j\in ds\right]\\[2mm] \notag
&& =:L_{11}(x_i,\,x_k,\,y_j)+L_{12}(x_i,\,x_k,\,y_j)+L_{21}(x_i,\,x_k,\,y_j)+L_{22}(x_i,\,x_k,\,y_j)\,, 
\eeam
Hence, by Assumption \ref{ass.KP.A}, we find
\beam \label{eq.KLP.5.3} \notag
&&L_{22}(x_i,\,x_k,\,y_j) \\[2mm] 
&& =\int_{\widehat{b}(x_i \wedge x_k)}^{\infty} \int_{c_j(y_j)}^{\infty}  \PP\left[|X_i|>\dfrac{x_i}t,\; X_k>\dfrac {x_k}t,\;Y_j>\dfrac {y_j}s\right]\, \PP\left[\widehat{\Theta}_{ik}\in dt,\;\Delta_j\in ds\right] \leq\\[2mm] \notag
&&\PP\left[\widehat{\Theta}_{ik}>\widehat{b}(x_i \wedge x_k),\;\Delta_j>c_j(y_j)\right]\leq\PP\left[\Delta_j>c_j(y_j)\right]=  o\left(\PP[\Theta_k\,X_k >x_k,\,\Delta_j\,Y_j>y_j]\right)\,,
\eeam
as $x_i\wedge x_k \wedge y_j \to \infty$. Next we calculate
\beam \label{eq.KLP.5.4} \notag
&&L_{21}(x_i,\,x_k,\,y_j) \\[2mm] \notag
&& =\int_{\widehat{b}(x_i \wedge x_k)}^{\infty} \int_{0}^{c_j(y_j)}  \PP\left[|X_i|>\dfrac{x_i}t,\; X_k>\dfrac {x_k}t,\;Y_j>\dfrac {y_j}s\right] \PP\left[\widehat{\Theta}_{ik}\in dt,\;\Delta_j\in ds\right]\\[2mm] 
&& \leq \PP\left[\widehat{\Theta}>\widehat{b}(x_i \wedge x_k),\;\Delta_j\leq c_j(y_j)\right]\leq \PP\left[\widehat{\Theta}_{ik}>\widehat{b}(x_i \wedge x_k)\right]\\[2mm] \notag
&&=o\left(\PP[\Theta_k\,X_k >x_k,\,\Delta_j\,Y_j>y_j]\right)\,,
\eeam
as $x_i\wedge x_k \wedge y_j \to \infty$, where in the last step we used relation \eqref{eq.KLP.5.1} . Similarly to $L_{22}(x_i,\,x_k,\,y_j)$, we obtain
\beam \label{eq.KLP.5.5} \notag
L_{12}(x_i,\,x_k,\,y_j)  &\leq& \PP\left[\widehat{\Theta}\leq \widehat{b}(x_i \wedge x_k),\;\Delta_j> c_j(y_j)\right]\leq \PP\left[\Delta_j>c_j(y_j)\right]\\[2mm] 
&=&o\left(\PP[\Theta_k\,X_k >x_k,\,\Delta_j\,Y_j>y_j]\right)\,,
\eeam
as $x_i\wedge x_k \wedge y_j \to \infty$. For the first term from the $GTAI$ property of $X_1,\,\ldots,\,X_n,\,Y_1,\,\ldots,\,Y_m$ and the definitions of $\widehat{b}$ and $c_j$ we have that
\beam \label{eq.KLP.5.6} \notag
L_{11}(x_i,\,x_k,\,y_j)  &\leq& \PP\left[|X_i|>\dfrac{x_i}{\widehat{b}(x_i \wedge x_k)},\; X_k>\dfrac {x_k}{\widehat{b}(x_i \wedge x_k)},\;Y_j>\dfrac {y_j}{c_j(y_j)}\right] \\[2mm] 
&=&o\left(\PP[X_k >x_k,\,Y_j>y_j]\right)\,,
\eeam
as $x_i\wedge x_k \wedge y_j \to \infty$, where in the last step we use that $(X_k,\,Y_j) \in \mathcal{D}^{(2)}$, for any $i=1,\,\ldots,\,n$, $j=1,\,\ldots,\,m$. Therefore by \eqref{eq.KLP.5.2}, together with \eqref{eq.KLP.5.3}, \eqref{eq.KLP.5.4}, \eqref{eq.KLP.5.5}, \eqref{eq.KLP.5.6}  we conclude that
\beam \label{eq.KLP.5.7}  \notag
&&\PP\left[|\Theta_i\,X_i|>x_i,\;\Theta_k\,X_k>x_k,\; \Delta_j\,Y_j>y_j\right]=o\left(\PP[\Theta_k\,X_k >x_k,\,\Delta_j\,Y_j>y_j]\right)  \\[2mm]
&&\qquad \qquad +o\left(\PP[X_k >x_k,\,Y_j>y_j]\right)=o\left(\PP[\Theta_k\,X_k >x_k,\,\Delta_j\,Y_j>y_j]\right)\,,
\eeam
as $x_i\wedge x_k \wedge y_j \to \infty$, where in last step we used that form assumptions follows
\beam \label{eq.KLP.5.8} 
\PP[\Theta_k\,X_k >x_k,\,\Delta_j\,Y_j>y_j] \asymp \PP[X_k >x_k,\,Y_j>y_j]\,,
\eeam
as $ x_k \wedge y_j \to \infty$. Indeed, since
\beam \label{eq.KLP.5.9} 
&&\PP[\Theta_k\,X_k >x_k,\,\Delta_j\,Y_j>y_j] \\[2mm] \notag
&&= \left(\int_0^{b_{k}(x_k)}+\int_{b_{k}(x_k)}^{\infty}\right) \left(\int_0^{c_j(y_j)}+\int_{c_j(y_j)}^{\infty}\right)  \PP\left[ X_k>\dfrac {x_k}t,\;Y_j>\dfrac {y_j}s\right] \PP\left[\Theta_{k}\in dt,\;\Delta_j\in ds\right]\\[2mm] \notag
&& =:M_{11}(x_k,\;y_j)+M_{12}(x_k,\;y_j)+M_{21}(x_k,\;y_j)+M_{22}(x_k,\;y_j)\,,
\eeam
From definitions of functions $b_{k}$ and $c_j$, see Assumption  \ref{ass.KP.A} and by $(X_k,\,Y_j) \in \mathcal{D}^{(2)}$ we find
\beam \label{eq.KLP.5.10} 
M_{11}(x_k,\;y_j)  \leq \PP\left[X_k >\dfrac{x_k}{b_k(x_k)},\,Y_j>\dfrac {y_j}{c_j(y_j)}\right]\,.
\eeam
Further we observe that
\beam \label{eq.KLP.5.11} 
M_{22}(x_k,\;y_j)&\leq& \PP\left[\Theta_{k} > b_{k}(x_k) ,\;\Delta_j>c_j(y_j)\right]\leq  \PP\left[\Delta_j>c_j(y_j)\right]   \\[2mm] \notag
&=&o(\PP[\Theta_k\,X_k >x_k,\,\Delta_j\,Y_j>y_j])\,,
\eeam
as $x_k \wedge y_j \to \infty$, where in the last step we employed Assumption  \ref{ass.KP.A}. Next, we see that 
\beam \label{eq.KLP.5.12}
M_{21}(x_k,\;y_j)&\leq& \PP\left[\Theta_{k} > b_{k}(x_k),\;\Delta_j\leq c_j(y_j)\right]\leq  \PP\left[\Theta_{k} > b_{k}(x_k) \right]  \\[2mm] \notag
& =&o(\PP[\Theta_k\,X_k >x_k,\,\Delta_j\,Y_j>y_j])\,,
\eeam
as $x_k \wedge y_j \to \infty$. Similarly by symmetry we obtain
\beam \label{eq.KLP.5.13} 
&&M_{12}(x_k,\;y_j) =o(\PP[\Theta_k\,X_k >x_k,\,\Delta_j\,Y_j>y_j])\,,
\eeam
as $x_k \wedge y_j \to \infty$. Hence, from relations \eqref{eq.KLP.5.10} -\eqref{eq.KLP.5.13}, together with \eqref{eq.KLP.5.9}, we get 
\beao
\PP[ \Theta_k\,X_k>x_k,\;\Delta_j\,Y_j>y_j] \lesssim \PP\left[ X_k>\dfrac {x_k}{b_k(x_k)},\;Y_j>\dfrac {y_j}{c_j(y_j)} \right] \asymp \PP[ X_k>x_k,\;Y_j>y_j]\,,
\eeao
as $ x_k \wedge y_j   \to \infty$. As a result we conclude that
\beam \label{eq.KP.5.13b}
\limsup_{x\wedge y \to \infty}\dfrac{\PP[ \Theta_k\,X_k>x_k,\;\Delta_j\,Y_j>y_j] }{\PP[X_k>x_k,\;Y_j>y_j] } < \infty\,.
\eeam

From the other side, for some arbitrarily chosen $\vep\in (0,\,1)$, we obtain
\beam \label{eq.KLP.5.1*} \notag 
&&\PP[\Theta_k\,X_k >x_k,\,\Delta_j\,Y_j>y_j]\\[2mm] \notag 
&&\geq \left(\int_{\vep}^1 + \int_1^{\infty} \right)\,\left( \int_{\vep}^1 + \int_1^{\infty} \right) \PP\left[X_k >\dfrac {x_k}{t},\,Y_j>\dfrac {y_j}{s}\right]\,\PP[\Theta_k \in dt\,,\;\Delta_j \in ds]\\[2mm] \notag 
&&\geq d_{1,1}\, \PP\left[X_k >x_k,\,Y_j>y_j \right]\,\PP[\Theta_k \in (\vep,\,1]\,,\;\Delta_j \in (\vep,\,1]]\\[2mm] \notag 
&&+\int_1^{\infty} \, \int_{\vep}^1 \PP\left[X_k >\dfrac {x_k}{t},\,Y_j>y_j \right]\,\PP[\Theta_k \in dt\,,\;\Delta_j \in d s]\\[2mm] \notag 
&&+ \int_{\vep}^1\,\int_1^{\infty}\PP\left[X_k >x_k,\,Y_j>\dfrac {y_j}{s} \right]\,\PP[\Theta_k \in dt\,,\;\Delta_j \in d s]\\[2mm] \notag 
&&+ d_{2,2}\PP\left[X_k >x_k,\,Y_j>y_j \right]\,\PP[\Theta_k \in (1,\,\infty)\,,\;\Delta_j \in (1,\,\infty)]\\[2mm] 
&&\geq d_{1,1} \PP\left[X_k >x_k,\,Y_j>y_j \right]\,\PP[\Theta_k \in (\vep,\,1]\,,\;\Delta_j \in (\vep,\,1]]\\[2mm] \notag 
&&+ d_{1,2} \PP\left[X_k >x_k,\,Y_j>y_j \right]\,\PP[\Theta_k \in (\vep,\,1]\,,\;\Delta_j \in (1,\,\infty)]\\[2mm] \notag 
&&+ d_{2,1} \PP\left[X_k >x_k,\,Y_j>y_j \right]\,\PP[\Theta_k \in (1,\,\infty)\,,\;\Delta_j \in (\vep,\,1]]\\[2mm] \notag 
&&+ d_{2,2} \PP\left[X_k >x_k,\,Y_j>y_j \right]\,\PP[\Theta_k \in (1,\,\infty)\,,\;\Delta_j \in (1,\,\infty)]\\[2mm] \notag 
&&\geq (d_{1,1}+d_{1,2}+d_{2,1}+d_{2,2})\, \PP\left[X_k >x_k,\,Y_j>y_j \right]\, \Big(\PP[\Theta_k \in (\vep,\,1]\,,\;\Delta_j \in (\vep,\,1]]\\[2mm] \notag 
&&+\PP[\Theta_k \in (\vep,\,1]\,,\;\Delta_j \in (1,\,\infty)]+ \PP[\Theta_k \in (1,\,\infty)\,,\;\Delta_j \in (\vep,\,1]]\\[2mm] \notag 
&&+\PP[\Theta_k \in (1,\,\infty)\,,\;\Delta_j \in (1,\,\infty)]\Big) \rightarrow (d_{1,1}+d_{1,2}+d_{2,1}+d_{2,2})\, \PP\left[X_k >x_k,\,Y_j>y_j \right]\,, 
\eeam
as $\vep \downarrow 0$, where the $d_{2,2}>0$ follows by class $\mathcal{D}^{(2)}$ property and the inclusion $(t,\,s) \in (1,\,\infty)\times (1,\,\infty)$. The inequalities $d_{1,1},\,d_{1,2},\,d_{2,1}>0$ follow from the intervals where the pairs $(t,\,s)$ belong. Hence by \eqref{eq.KLP.5.1*} we obtain
\beam \label{eq.KLP.5.2*}
\limsup_{x_k\wedge y_j \to \infty}\dfrac{\PP[ X_k>x_k,\;Y_j>y_j] }{\PP[\Theta_k\,X_k>x_k,\;\Delta_j\,Y_j>y_j] } < \infty\,,
\eeam
therefore, from relations \eqref{eq.KP.5.13b} and \eqref{eq.KLP.5.2*} we find \eqref{eq.KLP.5.8}.

With similar handling, because of symmetry, we obtain
\beam \label{eq.KLP.5.14}  
\PP[|\Delta_j\,Y_j|>y_j,\;\Delta_k\,Y_k>y_k,\; \Theta_i\,X_i>x_i]=o(\PP[\Theta_i\,X_i >x_i,\,\Delta_k\,Y_k>y_k])\,,
\eeam
as $x_i\wedge y_k \wedge y_j \to \infty$ for any $j\neq k=1,\,\ldots,\,m$, $i=1,\,\ldots,\,n$. So from relations \eqref{eq.KLP.5.7} and   \eqref{eq.KLP.5.14}, we find that the $\Theta_1\,X_1,\,\ldots,\,\Theta_n\,X_n$, and $\Delta_1\,Y_1,\,\ldots,\,\Delta_m\,Y_m$ are $GTAI$.

For the second case, with $GQAI$ structure, we follow the same route, with the only difference that the convergences are with $x \wedge y \to \infty$, and in relations  \eqref{eq.KLP.5.6} and \eqref{eq.KLP.5.7} the last term takes the form 
\beao
o(\PP[\Theta_i\,X_i >x,\;\Delta_j\,Y_j>y] + \PP[\Theta_k\,X_k >x,\;\Delta_j\,Y_j>y])\,.~\halmos
\eeao

Now we are ready to present the weighted form of Theorem \ref{th.KLP.3.1}, and the same time to generalize \cite[Th. 6.1]{konstantinides:passalidis:2025}, since the random weights $\Theta_i$, $\Delta_j$ are no more strictly positive and bounded from above and the main components $X_i$, $Y_j$ are real valued, random variables.

\bth \label{th.KLP.5.1}
\begin{enumerate}
\item[(i)]
Let $\{X_1,\,\ldots,\,X_n\}$, $\{Y_1,\,\ldots,\,Y_m\}$ be $GQAI$, real random variables with corresponding distributions $F_1,\,\ldots,\,F_n,\,G_1,\,\ldots,\,G_m$,  and for any $i=1,\,\ldots,\,n$, $j=1,\,\ldots,\,m$, it holds $(X_i,\,Y_j) \in \mathcal{C}^{(2)}$. We also assume that $\Theta_1,\,\ldots,\,\Theta_n$, $\Delta_1,\,\ldots,\,\Delta_m$ be non-negative,non-degenerate to zero random variables, independent of the main variables $X_1,\,\ldots,\,X_n,\,Y_1,\,\ldots,\,Y_m$, with $(\Theta_i,\,\Delta_j)$ satisfying Assumption \ref{ass.KP.A}. Then it holds
\beam \label{eq.KLP.5.15} \notag
\PP[S_n^{\Theta}> x\,,\;T_m^{\Delta}> y] \sim \PP\left[\bigvee_{i=1}^n S_i^{\Theta}> x\,,\; \bigvee_{j=1}^m T_j^{\Delta}> y\right] 
\sim  \sum_{i=1}^n \sum_{j=1}^m \PP[\Theta_i\,X_i>x\,,\;\Delta_j\,Y_j>y ]\,,\\
\eeam
as $x \wedge y \to \infty$.
\item[(ii)]
Let the assumptions of part $(1)$ hold, with the only differences that $n=m$, the $\{X_1,\,\ldots,\,X_n\}$, $\{Y_1,\,\ldots,\,Y_n\}$ are $GTAI$, real, random variables with corresponding distributions $F_1,\,\ldots,\,F_n,\,G_1,\,\ldots,\,G_n$, and for any $i,\,j=1,\,\ldots,\,n$,  hold $(X_i,\,Y_j) \in (\mathcal{D} \cap \mathcal{L})^{(2)}$. Further we assume that $\{X_1,\,\ldots,\,X_n\}$ are TAI and $\{Y_1,\,\ldots,\,Y_n\}$ are TAI and for any $i=1,\,\ldots,\,n$ $E[\Theta_{i}^{p}]<\infty$, for some $p>\bigvee_{i=1}^{n}J_{F_i}^+$, and for any $j=1,\,\ldots,\, n$ $E[\Delta_{j}^{q}]<\infty$, for some $q>\bigvee_{j=1}^{n}J_{G_j}^+$. Then it holds
\beam \label{eq.KLP.5.16} \notag
\PP\left[S_n^{\Theta}> x\,,\;T_n^{\Delta}> y\right] &\sim& \PP\left[\bigvee_{i=1}^n S_i^{\Theta}> x\,,\; \bigvee_{j=1}^n T_j^{\Delta}> y\right] \sim \PP\left[\bigvee_{i=1}^n \Theta_i\,X_{i}> x\,,\; \bigvee_{j=1}^n \Delta_j\,Y_{j}> y\right]\\[2mm] 
&\sim& \sum_{i=1}^n \sum_{j=1}^n \PP[\Theta_i\,X_i>x\,,\;\Delta_j\,Y_j>y ]\,,
\eeam
as $x \wedge y \to \infty$.
\end{enumerate}
\ethe

\pr~
\begin{enumerate}
\item[(i)]
Rewriting Theorem \ref{eq.KLP.3.1}, we find that $Z_1,\,\ldots,\,Z_n,\,W_1,\,\ldots,\,W_m$ are $GQAI$ with distributions from class $\mathcal{C}$ and for any $i=1,\,\ldots,\,n$, $j=1,\,\ldots,\,m$ hold $(Z_i,\,W_j) \in \mathcal{C}^{(2)}$, hence
\beam \label{eq.KLP.5.17} \notag
\PP\left[\sum_{i=1}^n Z_{i}> x\,,\;\sum_{j=1}^m W_{j}> y\right] &\sim& \PP\left[\bigvee_{k=1}^n \sum_{i=1}^k Z_{i}> x\,,\; \bigvee_{l=1}^m \sum_{j=1}^l W_{j}> y\right] \\[2mm]
&\sim&  \sum_{i=1}^n \sum_{j=1}^m \PP[Z_i>x\,,\;W_j>y ]\,,
\eeam
as $x \wedge y \to \infty$. From the assumptions of part $(i)$, applying Theorem \ref{th.KLP.4.2}(iii), since for any $i=1,\,\ldots,\,n$, $j=1,\,\ldots,\,m$ holds $(X_i,\,Y_j) \in \mathcal{C}^{(2)}$, it follows that $(\Theta_i\,X_i,\,\Delta_j\,Y_j) \in \mathcal{C}^{(2)}$. Further, applying Lemma \ref{lem.KLP.5.1} we find out that the products $\Theta_1\,X_1,\,\ldots,\,\Theta_n\,X_n$,  $\Delta_1\,Y_1,\,\ldots,\,\Delta_m\,Y_m$ are $GQAI$, as $\mathcal{C}^{(2)} \subsetneq \mathcal{D}^{(2)}$. Hence, putting 
\beam \label{eq.K.5.18}
Z_i:= \Theta_i\,X_i\,, \qquad W_j:= \Delta_j\,Y_j\,,
\eeam 
for any $i=1,\,\ldots,\,n$ and any $j=1,\,\ldots,\,m$, we find \eqref{eq.KLP.5.15} through \eqref{eq.KLP.5.17}. 
\item[(ii)]
Repeating relation \eqref{eq.KLP.1.8} through the $Z_1,\,\ldots,\,Z_n,\,W_1,\,\ldots,\,W_n$ that are $GTAI$ (and each one of these sequences are TAI)  with any $i,\,j=1,\,\ldots,\,n$, because of $(Z_i,\,W_j) \in (\mathcal{D} \cap \mathcal{L})^{(2)}$ we obtain
\beam \label{eq.KLP.5.18}
&&\PP\left[\sum_{i=1}^n Z_{i}> x\,,\;\sum_{j=1}^n W_{j}> y\right] \sim \PP\left[\bigvee_{k=1}^n \sum_{i=1}^k Z_{i}> x\,,\; \bigvee_{l=1}^n \sum_{j=1}^l W_{j}> y\right] \\[2mm] \notag
&&\qquad \qquad \qquad \qquad \sim\PP\left[\bigvee_{i=1}^n Z_{i}> x\,,\; \bigvee_{j=1}^n W_{j}> y\right] \sim \sum_{i=1}^n \sum_{j=1}^m \PP[Z_i>x\,,\;W_j>y ]\,,
\eeam
as $x \wedge y \to \infty$, see also \cite[Th. 4.2]{konstantinides:passalidis:2024b}. Therefore, since $(X_i,\,Y_j) \in (\mathcal{D} \cap \mathcal{L})^{(2)}$, from Theorem \ref{th.KLP.4.2}(iv), because of Assumption \ref{ass.KP.A}, we conclude  $(\Theta_i\,X_i,\,\Delta_j\,Y_j) \in (\mathcal{D} \cap \mathcal{L})^{(2)}$. Further by Lemma \ref{lem.KLP.5.1} we obtain that the $\Theta_1\,X_1,\,\ldots,\,\Theta_n\,X_n,\,\Delta_1\,Y_1,\,\ldots,\,\Delta_n\,Y_n$ are $GTAI$, and by the moment condition and the class $\mathcal{D}$ (for the primary random variables) we have that $\Theta_1\,X_1,\,\ldots,\,\Theta_n\,X_n$ are TAI and $\Delta_1\,Y_1,\,\ldots,\,\Delta_n\,Y_n$ are TAI, by Theorem 2.2 of \cite{li:2013}. Thus, using \eqref{eq.K.5.18} for any $i,\,j=1,\,\ldots,\,n$, we find \eqref{eq.KLP.5.16} through \eqref{eq.KLP.5.18}.  ~\halmos
\end{enumerate}

\bre \label{rem.KLP.5.1}
Comparing the two parts of Theorem \ref{th.KLP.5.1} we realize that as the distribution class $\mathcal{C}^{(2)}$ increases to $(\mathcal{D} \cap \mathcal{L})^{(2)}$, the dependence decreases from $GQAI$ to $GTAI$ and additionally we need $n=m$. However, it is remarkable, that in $GTAI$ case we find the asymptotic behavior of the maximums, which is NOT possible in the part $(i)$.
Furthermore in the second part we obtain
\beam \label{eq.KLP.5.19}
\PP\left[\bigvee_{i=1}^n \Theta_{i}\,X_i> x\,,\; \bigvee_{j=1}^m \Delta_{j}\,Y_j> y\right] \sim \sum_{i=1}^n \sum_{j=1}^m \PP[\Theta_i\,X_i>x\,,\;\Delta_j\,Y_j>y ]\,,
\eeam
as $x \wedge y \to \infty$, whose proof follows the same arguments of Theorem \ref{th.KLP.5.1}(ii), (and without necessary TAI for each sequence!) namely we establish \eqref{eq.KLP.5.19} via closure property of $(\mathcal{D} \cap \mathcal{L})^{(2)}$ and with respect to $GTAI$ structure of the products, in combination with \cite[Th. 4.1]{konstantinides:passalidis:2025}.
\ere

The next result, with an extra condition on the moments of the random weights, provides a direct formula for calculation in the distribution class $\mathcal{R}_{(-\alpha_1,\,-\alpha_2)}^{(2)}$.

\bco \label{cor.KLP.5.1}
\begin{enumerate}
\item[(i)]
Let hold the assumptions of Theorem \ref{th.KLP.5.1}(i), with the only difference $F_i \in \mathcal{R}_{-\alpha_{1i}}$, $G_j \in \mathcal{R}_{-\alpha_{2j}}$ and $(X_i,\,Y_j) \in \mathcal{R}_{(-\alpha_{1i},\,-\alpha_{2j})}^{(2)}$ with $0<\alpha_{1i},\,\alpha_{2j} < \infty$ for any $i=1,\,\ldots,\,n$, $j=1,\,\ldots,\,m$. Additionally we assume that 
\beao
\E\left[\Theta_i^{\alpha_{1i}+\vep}\,\Delta_j^{\alpha_{2j}+\vep}\right]< \infty\,,
\eeao 
for any $i=1,\,\ldots,\,n$ and $j=1,\,\ldots,\,m$, for some $\vep>0$. Then
\beao
\PP[S_n^{\Theta}> x\,,\;T_m^{\Delta}> y] &\sim& \PP\left[\bigvee_{i=1}^n S_i^{\Theta}> x\,,\; \bigvee_{j=1}^m T_j^{\Delta}> y\right] \\[2mm]
&\sim&  \sum_{i=1}^n \sum_{j=1}^m \E\left[\Theta_i^{\alpha_{1i}}\,\Delta_j^{\alpha_{2j}}\right]\,\PP[X_i>x\,,\;Y_j>y ]\,.
\eeao
as $x \wedge y \to \infty$.
\item[(ii)]
Let hold the assumptions of Theorem \ref{th.KLP.5.1}(ii), with the only difference  $F_i \in \mathcal{R}_{-\alpha_{1i}}$, $G_j \in \mathcal{R}_{-\alpha_{2j}}$ and $(X_i,\,Y_j) \in \mathcal{R}_{(-\alpha_{1i},\,-\alpha_{2j})}^{(2)}$with $0<\alpha_{1i},\,\alpha_{2j} < \infty$, for any $i,\,j=1,\,\ldots,\,n$. Furthermore we assume that $\E\left[\Theta_i^{\alpha_{1i}+\vep}\,\Delta_j^{\alpha_{2j}+\vep}\right]< \infty$, for some $\vep>0$, for any $i,\,j=1,\,\ldots,\,n$. Then
\beao
\PP[S_n^{\Theta}> x\,,\;T_n^{\Delta}> y] &\sim& \PP\left[\bigvee_{i=1}^n S_i^{\Theta}> x\,,\; \bigvee_{j=1}^n T_j^{\Delta}> y\right] \sim \PP\left[\bigvee_{i=1}^n \Theta_i\,X_{i}> x\,,\; \bigvee_{j=1}^n \Delta_j\,Y_{j}> y\right]\\[2mm] 
&\sim& \sum_{i=1}^n \sum_{j=1}^n \E\left[\Theta_i^{\alpha_{1i}}\,\Delta_j^{\alpha_{2j}}\right]\,\PP[X_i>x\,,\;Y_j>y ]\,,
\eeao
as $x \wedge y \to \infty$.
\end{enumerate}
\eco

\pr~
The arguments follow directly from application of Lemma \ref{lem.KLP.4.2} on Theorem \ref{th.KLP.5.1}, since $\mathcal{R}_{(-\alpha_{1i},\,-\alpha_{2j})}^{(2)} \subsetneq \mathcal{C}^{(2)} \subsetneq (\mathcal{D}\cap \mathcal{L})^{(2)}$. ~\halmos

Now we apply these results to two-dimensional discrete time risk model with stochastic returns. The surplus processes take the form 
\beao
U_1^{\Theta}(k,\,x):=x - \sum_{i=1}^k \Theta_i\,X_i\,, \qquad U_2^{\Delta}(k,\,y):=y - \sum_{j=1}^k \Delta_j\,Y_j\,, 
\eeao
for some $k=1,\,\ldots,\,n$, where $x$ and $y$ are the initial surpluses, in each business lines, the $X_i$, $Y_j$ are the net losses in the $i$-th and $j$-th period respectively, and the $\Theta_i$ and $\Delta_j$ play the role of discount factor. 

Let choose the following definitions \texttt{}
\beao
\tau_1^{\Theta}(x):= \inf \{k=1,\,\ldots,\,n\;:\; U_1^{\Theta}(k,\,x)<0 \;|\; U_1^{\Theta}(0)=x\}\,,\\[2mm]
\tau_2^{\Delta}(y):= \inf \{k=1,\,\ldots,\,n\;:\; U_2^{\Delta}(k,\,y)<0 \;|\; U_2^{\Delta}(0)=y\}\,,
\eeao
and $\tau_{and}^{\Theta,\Delta}:=\max \left\{\tau_1^{\Theta}(x),\,\tau_2^{\Delta}(y)\right\}$, that indicates the first moment, when both surpluses fall below zero, but not necessarily simultaneously. With these notations we provide the formula for the ruin probability in two-dimensional discrete risk model, over finite time horizon
\beam \label{eq.KLP.5.23} 
\psi_{and}^{\Theta,\Delta}(x,\,y;\,n):=\PP\left[\tau_{and}^{\Theta,\Delta}\ \leq n\right]=\PP\left[\bigvee_{i=1}^n S_i^{\Theta}> x\,,\; \bigvee_{j=1}^n T_j^{\Delta}> y\right]\,.
\eeam   

\bco \label{cor.KLP.5.2}
\begin{enumerate}
\item[(i)]
Let hold the assumptions of Theorem \ref{th.KLP.5.1}[(i) or (ii)]. Then
\beao
\psi_{and}^{\Theta,\Delta}(x,\,y;\,n) \sim  \sum_{i=1}^n \sum_{j=1}^n \PP[\Theta_i\,X_i>x\,,\;\Delta_j\,Y_j>y ]\,.
\eeao
as $x \wedge y \to \infty$.
\item[(ii)]
Let hold the assumptions of Corollary \ref{cor.KLP.5.1}[(i) or (ii)]. Then
\beao
\psi_{and}^{\Theta,\Delta}(x,\,y;\,n)  \sim \sum_{i=1}^n \sum_{j=1}^n \E\left[\Theta_i^{\alpha_{1i}}\,\Delta_j^{\alpha_{2j}}\right]\,\PP[X_i>x\,,\;Y_j>y ]\,,
\eeao
as $x \wedge y \to \infty$.
\end{enumerate}
\eco

\pr~
By the results from (i) and (ii) are implied from Theorem \ref{th.KLP.5.1} and Corollary \ref{cor.KLP.5.1} respectively we get the desired concluson.
~\halmos

\noindent \textbf{Acknowledgments.} 
We feel the pleasant duty to express deep gratitude to prof. Jinzhu Li for his useful comments that improved significantly the text.


\begin{thebibliography}{99}




\bibitem{asmussen:2003}
{\sc Asmussen, S.}\ (2003)
{\em Applied Probability and Queues} 
Springer, New York, 2nd ed.


\bibitem{asmussen:klueppelberg:1996}
{\sc Asmussen, S., Kl\"{u}ppelberg, C.}\ (1996)
Large deviations results for subexponential tails, with applications to insurance risk. 
{\em Stoch. Process. Appl.}, \textbf{64}, 103--125.



\bibitem{athreya:ney:1972}
{\sc Athreya, K.B., Ney, P.E.}\ (1972)
{\em Branching Processes} 
Springer, New York.



\bibitem{basrak:davis:mikosch:2002}
{\sc Basrak, B., Davis, R.A., Mikosch, T.}\ (2002) 
A characterization of multivariate regular variation.
{\em  Ann. Appl. Probab.} \textbf{12}, 908--920.

\bibitem{basrak:davis:mikosch:2002b}
{\sc Basrak, B., Davis, R.A., Mikosch, T.}\ (2002) 
Regular variation of GARCH processes.
{\em  Stoch. Process. Appl.} \textbf{99}, no. 1, 95--115.



\bibitem{breiman:1965}
{\sc Breiman L.}\ (1965) 
On some limit theorems similar to arc-sin law.
{\em Theory Probab. Appl.}, \textbf{10}, 323--331.

\bibitem{buraczewski:damek:mikosch:2016}
{\sc Buraczewski, D., Damek, E., Mikosch, T.}\ (2016)
{\em Stochastic Models with Power-Law Tails} 
Springer, New York.


\bibitem{cai:tang:2004}
{\sc Cai, J., Tang, Q.}\  (2004)
On max-sum equivalence and convolution closure of heavy-tailed distributions and their applications. 
{\em J. Appl. Probab.} \textbf{41}, 117--130.

\bibitem{chen:cheng:zheng:2025}
{ \sc Chen, Z., Cheng, D., Zheng, H.}\ (2025)
On the joint tail behavior of randomly weighted sums of dependent random variables with applications to risk theory.
{ \em Scand. Acutar. J.}, 1--20.






\bibitem{chen:yang:2019}
{\sc Chen, Y., Yang, Y.}\  (2019)
Bi-variate regular variation among randomly weighted sums in general insurance. 
{\em Eur. Actuar. J.}, \textbf{9}, 301--322.

\bibitem{chen:yuen:2009}
{ \sc Chen, Y., Yuen, K.C.}\ (2009)
Sums of pairwise quasi-asymptotic independent random variables with consistent variation.
{ \em Stochastic Models}, \textbf{25}, 76--89.

\bibitem{chen:wang:wang:2013}
{ \sc Chen, Y., Wang, L., Wang, Y.}\ (2013)
Uniform asymptotics for the finite-time ruin probabilities of two kinds of nonstandard bidimensional risk models.
{ \em J. Math. Anal. Appl.}, \textbf{401}, no. 1, 114--129.



\bibitem{cheng:yu:2019}
{\sc Cheng, D., Yu, C.}\  (2019)
Uniform asymptotics for the ruin probabilities in a bidimensional renewal risk model with strongly subexponential claims.
{\em Stochastics} \textbf{91}, Vol 1. 643--656.


\bibitem{cheng:konstantinides:wang:2024}
{ \sc Cheng, M., Konstantinides, D.G., Wang, D}\ (2024)
Multivariate regular varying insurance and financial risks in $d$-dimensional risk model. 
{ \em J. Appl. Probab.},  \textbf{61}, no. 4, 1319 -- 1342.

\bibitem{chistyakov:1964}
{ \sc Chistyakov, V.P.}\ (1964)
A theorem on sums of independent positive random variables and its applications to branching random processes.
{ \em Theory Probab. Appl.}, \textbf{9}, 640--648.

\bibitem{cline:resnick:1992} 
{\sc Cline, D.B.H., Resnick, S.}\ (1992)
Multivariate subexponential distributions.
{\em Stoch. Process. Appl.}, \textbf{42}, no.1, 49--72.

\bibitem{cline:samorodnitsky:1994} 
{\sc Cline, D.B.H., Samorodnitsky, G.}\ (1994)
Subexponentiality of the product of independent random variables.
{\em Stoch. Process. Appl.}, \textbf{49}, 75--98.

\bibitem{cui:wang:2020} 
{\sc Cui, Z., Wang, Y.}\ (2020)
On the long tail property of product convolution.
{\em Lith. Math. J.}, \textbf{60}, no. 2, 315--329.
\bibitem{das:fasenhartmann:2023}
{\sc Das, B., Fasen-Hartmann, V.}\ (2023)
Aggregating heavy-tailed random vectors: from finite sums to Levy processes.
 {\em arXiv:2301.10423v1}.

\bibitem{denisov:zwart:2007}
{\sc Denisov, D., Zwart, B.}\ (2007)
 On a theorem of Breiman and a class of random difference equations.
 {\em J. Appl. Probab.} \textbf{44}, 1031--1046.

\bibitem{dirma:nakiliuda:siaulys:2023} 
{\sc Dirma, M., Nakliuda, N., {\v{S}}iaulys, J.}\ (2023)
Generalized moments of sums with heavy-tailed random summands.
{\em  Lith. Math. J.}, \textbf{63}, no. 3, 254--271.






\bibitem{feller:1969}
{\sc Feller, W.}\ (1969)
One-sided analogues of  Karamata's regular  variation.
{\em L' enseignement Math\'{e}matique}, \textbf{15}, 107--121.

\bibitem{feller:1971} 
Feller, W. (1971)
\emph{ An Introduction to Probability Theory and Its Applications} Vol. II, second edition, 
Wiley, New York.

\bibitem{foss:korshunov:palmowski:2024} 
{\sc Foss, S., Korshunov, D., Palmowski, Z.} \ (2024)
Maxima over random  time intervals for heavy-tailed compound renewal and L\'{e}vy processes.
{\em Stoch. Process. Appl.}, \textbf{176}, 104422.

\bibitem{foss:korshunov:zachary:2013} 
{\sc Foss, S., Korshunov, D., Zachary, S.} \ (2013)
{\em An Introduction to Heavy-Tailed and Subexponential Distributions.}
Springer, New York, 2nd ed.

\bibitem{fougeres:mercadier:2012}
{\sc Fougeres, A., Mercadier, C.}\ (2012)
Risk measures and multivariate extensions of Breiman’s theorem.
{\em J. Appl. Probab.},  \textbf{49}, no. 2, 364--384.



\bibitem{geluk:tang:2009}
{\sc Geluk, J., Tang, Q.}\ (2009)
Asymptotic tail probabilities of sums of dependent subexponential random variables. 
{\em J. Theor. Probab.}, \textbf{22}, 871--882.


\bibitem{goldie:1978}
{\sc Goldie, C.M.}\ (1978)
Subexponential distributions and dominated variation tails 
{\em J. Appl. Probab.}, \textbf{15}, 440--442.

\bibitem{dehaan:resnick:1982} 
{\sc Haan, L. de, Resnick, S.}\ (1981)
On the observation closet to the origin.
{\em  Stoch. Process. Appl.}, \textbf{11}, no. 3, 301--308.





\bibitem{jiang:wang:chen:xu:2015}
{\sc Jiang, T., Wang, Y., Chen, Y., Xu, H.}\  (2015)
Uniform asymptotic estimate for finite-time ruin probabilities of a time-dependent bidimensional renewal model. 
{\em Insur. Math. Econom.}, \textbf{64}, 45--53.


\bibitem{kizinevic:sprindys:siaulys:2016}
{\sc Kizinevi{\v{c}}, E., Sprindys, J., {\v{S}}iaulys, J.}\ (2016)
Randomly stopped sums with consistently varying distributions. 
{\em Modern Stochastics: Theory and Applications}, \textbf{3}, 165--179.
	
\bibitem{ko:tang:2008} 
{\sc Ko, B.W., Tang Q.H.} \ (2008)
Sums of dependent nonnegative random variables with subexponential tails.
{\em J. Appl. Probab.}, \textbf{45}, 85--94.

\bibitem{konstantinides:2018} 
{\sc Konstantinides, D.G.} \ (2018)
{\em Risk Theory. A Heavy Tail Approach.}
World Scientific, New Jersey.


\bibitem{konstantinides:leipus:siaulys:2023}
{\sc Konstantinides, D.G., Leipus, R., \v{S}iaulys, J.}\ (2023) 
On the non-closure property for strong subexponential distributions. 
{\em Nonlinear Analysis: Modelling and Control} \textbf{28}, no. 1, 97--115.

\bibitem{konstantinides:leipus:passalidis:siaulys:2025}
{\sc Konstantinides, D.G., Leipus, R., Passalidis, C.D., \v{S}iaulys, J.}\ (2025) 
Tail behavior of randomly weighted sums with interdependent summands. 
{\em Preprint, arXiv: 2503.11271}.

\bibitem{konstantinides:li:2016}
{\sc Konstantinides, D.G., Li, J.}\ (2016)
Asymptotic ruin probabilities for a multidimensional renewal risk model with multivariate regularly varying claims.
{\em Insur. Math. Econom.}, \textbf{69}, 38--44.

\bibitem{konstantinides:mikosch:2005}
{\sc Konstantinides, D.G., Mikosch, T.}\ (2005)
Large Deviations and Ruin Probabilities for Solutions to Stochastic Recurrence Equations with Heavy-tailed Innovations.
{\em Ann. Probab.}, \textbf{33}, 1.992--2.035.

\bibitem{konstantinides:passalidis:2023} 
{\sc Konstantinides, D.G., Passalidis, C.D.} \ (2024)
Background risk model in presence of heavy tails under dependence.
{\em Preprint, arXiv:2405.03014}.

\bibitem{konstantinides:passalidis:2024b} 
{\sc Konstantinides, D.G., Passalidis, C.D.} \ (2024)
Closure properties and heavy tails: random vectors in the presence of dependence
{\em Preprint, arXiv:2402.09041}.

\bibitem{konstantinides:passalidis:2024g} 
{\sc Konstantinides, D.G., Passalidis, C.D.} \ (2024)
Random vectors in the presence of a single big jump.
{\em Preprint, arXiv:2410.10292}.

\bibitem{konstantinides:passalidis:2025} 
{\sc Konstantinides, D.G., Passalidis, C.D.} \ (2025)
A new approach in two-dimensional heavy-tailed distributions.
{\em Ann. Actuar. Scienc.}, 1--33.


\bibitem{korshunov:2018}
{\sc Korshunov, D.}\ (2018)
On subexponential tails for the maxima of negatively driven compound renewal and L\'{e}vy processes.
{\em Stoch. Process. Appl.}, \textbf{128}, 1.316--1.332.


\bibitem{leipus:siaulys:2020}
{\sc Leipus, R.,  {\v S}iaulys, J.}\ (2020)
On a closure property of convolution equivalent class of distributions.
{\em J. Math. Anal. Appl.}, \textbf{490}, no. 124226.



\bibitem{leipus:siaulys:konstantinides:2023}
{\sc Leipus, R., \v{S}iaulys, J., Konstantinides, D.G.}\ (2023)
{\em Closure Properties for Heavy-Tailed and Related Distributions: An Overview.}
Springer Nature, Cham.


\bibitem{leipus:surgailis:2007}
{\sc Leipus, R., Surgailis, D.}\ (2007)
On long-range dependence in regenerative processes based on a general ON/OFF scheme. 
{\em J. Appl. Probab.}, \textbf{44}, 379--392. 

\bibitem{li:2013}
{\sc Li, J.}\ (2013)
On pairwise quasi-asymptotically independent random variables and  their applications. 
{\em Stat. Prob. Lett.}, \textbf{83}, 2081--2087.


\bibitem{li:2017}
{\sc Li, J.}\ (2017)
A note on the finite-time ruin probabilities of a renewal risk model with Brownian perturbation. 
{\em Stat. Prob. Lett.}, \textbf{127}, 49--55.

\bibitem{li:2018a}
{\sc Li, J.}\ (2018)
A revisit to asymptotic ruin probabilities of a two-dimesional renewal risk model. 
{\em Stat. Probab. Lett.}, \textbf{140}, 23--32.

\bibitem{li:2018b}
{\sc Li, J.}\ (2018b)
On the joint tail behavior of randomly weighted sums of heavy-tailed random variables.
{\em J. Mutlivar. Anal.}, \textbf{164}, 40--53.


\bibitem{li:yang:2015}
{\sc Li, J., Yang, H.}\ (2015)
Asymptotic ruin probabilities for a bidimensional renewal risk model with constant interest rate and dependent claims.
{\em J. Math. Anal. Appl.}, \textbf{426}, no. 1, 247--266.





\bibitem{matuszewska:1964}
{\sc Matuszewska, W.}\ (1964)
On generalization of regularly increasing functions.
{\em Studia Mathematica}, \textbf{24}, 271--279.

\bibitem{maulik:resnick:2004} 
{\sc Maulik, K., Resnick, S.I.}\ (2004)
Characterizations and examples of hidden regular variation.
{\em Extremes}, \textbf{7}, 31--67.



\bibitem{omey:2006}
{\sc Omey, E.} \ (2006)
Subexponential distribution functions in $R^{d}$.
{\em J. Math. Sci.}, \textbf{138}, no. 1, 5434--5449.



\bibitem{pakes:2004}
{\sc Pakes, A.G.} \ (2004)
Convolution equivalence and infinite divisibility.
{\em J. Appl. Probab.}, \textbf{41}, 407--424.

\bibitem{resnick:2007}
{\sc Resnick, S.}\ (2007) 
{\em Heavy-Tail Phenomena. Probabilistic and Statistical Modeling.} 
Springer, New York.

\bibitem{rolski:schmidli:schmidt:teugels:1999}
{\sc Rolski, T., Schmidli, H., Schmidt, V., Teugels, J.L.}\ (1999) 
{\em Stochastic Processes for Insurance and Finance} 
Wiley, Chichester.



\bibitem{samorodnitsky:sun:2016} 
{\sc Samorodnitsky, G., Sun, J.}\ (2016) 
Multivariate subexponential distributions and their applications. 
{\em Extremes}, \textbf{19}, no. 2, 171--196.

\bibitem{samorodnitsky:taqqu:1994} 
{\sc Samorodnitsky, G., Taqqu, M.}\ (1994) 
{\em Stable Non-Gaussian Random Processes: Stochastic Models with Infinite Variance}
Chapman and Hall, New York.



\bibitem{shen:du:2023}
{\sc Shen, X., Du, K.}\ (2023) 
Uniform approximation for the tail behavior of bidimensional randomly weighted sums.
{\em Methodol. Comput. Appl. Probab.}, \textbf{25}, no. 26.

\bibitem{shen:ge:fu:2020}
{\sc Shen, X., Ge, M., Fu, K.A.}\ (2020) 
Approximation of the tail probabilities for bidimensional randomly weighted sums with dependent components.
{\em Probab. Eng. Inf. Sci.}, \textbf{34}, no. 1, 112--130.


\bibitem{tang:2006}
{\sc Tang, Q.}\ (2006)
The subexponentiality of products revisited.
{\em Extremes}, \textbf{9}, 231--241.

\bibitem{tang:2008}
{\sc Tang, Q.}\ (2008)
Insensitivity to negative dependence of asymptotic tail probabilities of sums and maxima of sums.
{\em Stoch. Anal. Appl.}, \textbf{26}, 435--450.



\bibitem{tang:tsitsiashvili:2003}
{\sc Tang, Q., Tsitsiashvili, G.}\ (2003)
Randomly weighted sums of subexponential random variables with application to ruin theory.
{\em Extremes}, \textbf{6}, 171--188.

\bibitem{tang:yuan:2014}
{\sc Tang, Q., Yuan, Z.}\ (2014) 
Randomly weighted sums of subexponential random variables with application to capital allocation.
{\em Extremes}, \textbf{17}, 467--493.

\bibitem{tang:yang:2019}
{\sc Tang, Q., Yang, Y.}\ (2019)
Interplay of insurance and financial risks in a stochastic environment. 
{\em Scand. Actuar. J.}, no. 5, 432--451.


\bibitem{wang:2011}
{\sc Wang, K.}\ (2011)
Randomly weighted sums of dependent subexponential random variables.
{\em Lith. Math. J.}, \textbf{51}, no. 4, 573--586.

\bibitem{wang:su:yang:2024}
{\sc Wang, H., Su, Q., Yang, Y.}\ (2024)
Asymptotics for ruin probabilities in a bidimensional discrete-time risk model with dependence and consistently varying net losses. 
{\em Stochastics}, \textbf{96}, no. 1, 667--695.

\bibitem{yang:chen:yuen:2024}
{\sc Yang, Y., Chen, S., Yuen, C.}\ (2024)
Asymptotics for the joint tail probability of bidimensional randomly weighted sums with applications to insurance.
{\em Sci. China Math.}, \textbf{67}, 163--186.



\bibitem{yang:leipus:siaulys:2012}
{\sc Yang, Y., Leipus, R., \v{S}iaulys, J.}\ (2012)
Tail probability of randomly weighted sums of subexponential random variables under a dependence structure. 
{\em Stat. Probab. Lett.}, \textbf{82}, 1727--1736.

\bibitem{yang:su:2023}
{\sc Yang, Y., Su, Q.}\ (2023)
Asymptotic behavior of ruin probabilities in a multidimensional risk model with investment and multivariate regularly varying claims.
{\em J. Math. Anal. Appl.}, \textbf{525}, 127319.



\bibitem{yang:wang:leipus:siaulys:2011}
{\sc Yang, Y., Wang, K., Leipus, R., \v{S}iaulys, J.}\ (2011)
Tail behavior of sums and maxima of sums of dependent subexponential random variables.
{\em Acta Appl. Math.}, \textbf{114}, 219--231.


\bibitem{yang:yuen:liu:2018}
{\sc Yang, Y., Yuen, K.C., Liu, J.-F.}\ (2018) 
Asymptotics for ruin probabilities in L\'{e}vy-driven risk models with heavy-tailed claims. 
{\em J. Ind. Manag. Optim.}, \textbf{14}, 231--247.

\bibitem{yuan:lu:2023}
{\sc Yuan, M., Lu, D.}\ (2023)
Asymptotics for a time-dependent by-claim model with dependent subexponential claims.
{\em Insur. Math. Econom.}, \textbf{112}, 120--141.


\bibitem{xu:shen:wang:2025}
{\sc Xu, C., Shen, X., Wang, K.}\ (2025)
The finite-time ruin probabilities of a dependent bidimensional risk model with subexponential claims and Brownian perturbations.
{\em Non. Anal. Mod. Contr.}, \textbf{30}, no. 3, 460--482.



\end{thebibliography}
\end{document}